\documentclass[12pt]{article}
\usepackage{latexsym,amssymb,amsmath}
\begin{document}

\baselineskip=20pt

\newcommand{\la}{\langle}
\newcommand{\ra}{\rangle}
\newcommand{\psp}{\vspace{0.4cm}}
\newcommand{\pse}{\vspace{0.2cm}}
\newcommand{\ptl}{\partial}
\newcommand{\dlt}{\delta}
\newcommand{\sgm}{\sigma}
\newcommand{\al}{\alpha}
\newcommand{\be}{\beta}
\newcommand{\G}{\Gamma}
\newcommand{\gm}{\gamma}
\newcommand{\vs}{\varsigma}
\newcommand{\Lmd}{\Lambda}
\newcommand{\lmd}{\lambda}
\newcommand{\td}{\tilde}
\newcommand{\vf}{\varphi}
\newcommand{\yt}{Y^{\nu}}
\newcommand{\wt}{\mbox{wt}\:}
\newcommand{\rd}{\mbox{Res}}
\newcommand{\for}{\mbox{for}}
\newcommand{\ad}{\mbox{\small ad}}
\newcommand{\voa}{(V,Y(\cdot,z),{\bf 1},\omega)}
\newcommand{\stl}{\stackrel}
\newcommand{\itp}{[^{\;\;\;W_3}_{\;W_1W_2}]}
\newcommand{\fr}{N^{\;\;W_3}_{W_1W_2}}
\newcommand{\ol}{\overline}
\newcommand{\ul}{\underline}
\newcommand{\es}{\epsilon}
\newcommand{\dmd}{\diamond}
\newcommand{\clt}{\clubsuit}
\newcommand{\vt}{\vartheta}
\newcommand{\ves}{\varepsilon}
\newcommand{\dg}{\dagger}
\newcommand{\tr}{\mbox{Tr}}
\newcommand{\ga}{{\cal G}({\cal A})}
\newcommand{\hga}{\hat{\cal G}({\cal A})}

\begin{center}{\LARGE \bf Quadratic Conformal Superalgebras}
\end{center}
\vspace{0.2cm}

\begin{center}{\large Xiaoping Xu}\end{center}
\begin{center}{Department of Mathematics, The Hong Kong University of Science \& Technology}\end{center}
\begin{center}{Clear Water Bay, Kowloon, Hong Kong}\footnote{Research supported
 by Hong Kong RGC Competitive Earmarked Research Grant HKUST709/96P.}\end{center}

\vspace{0.3cm}

\begin{center}{\Large \bf Abstract}\end{center}
\vspace{0.2cm}

{\small In this paper, we shall classify  ``quadratic'' conformal superalgebras by certain compatible pairs of a Lie superalgebra and a Novikov superalgebra. Four general constructions of such pairs  are given. Moreover, we shall classify  such pairs related to simple Novikov algebras.}

\section{Introduction}

The notion of conformal superalgebra was formulated by Kac [Ka2], which is equivalent to the notion of linear Hamiltonian operator in Gel'fand-Dikii-Dorfman's theory in [GDi1-2] and [GDo1-3]. Conformal superalgebras play important roles in quantum field theory (e.g. cf. [Ka2]) and vertex operator superalgebras (e.g. cf. [Ka2], [X7]). In some sense, conformal superalgebras are generalizations of affine Kac-Moody algebras and the Virasoro algebra. In this paper, we shall study a special class of conformal superalgebras, which we call ``quadratic conformal superalgebras.'' Below, we give a more detailed introduction.

Throughout this paper, all the vector spaces are assumed over $\Bbb{C}$, the field of complex numbers. Denote by $\Bbb{C}^+$ the additive group of $\Bbb{C}$. For two vector spaces $V_1$ and $V_2$, we denote by $LM(V_1,V_2)$ the space of linear maps from $V_1$ to $V_2$. Moreover, we denote by $\Bbb{Z}$ the ring of integers and by $\Bbb{Z}_2=\Bbb{Z}/2\Bbb{Z}$ the cyclic group of order 2. When the context is clear, we use $\{0,1\}$ to denote the elements of $\Bbb{Z}_2$. We shall also use the following operator of taking residue:
$$\rd_z(z^n)=\dlt_{n,-1}\qquad\for\;\;n\in \Bbb{Z}.\eqno(1.1)$$
Furthermore, all the binomials are assumed to be expanded in the nonnegative powers of the second variable. 

 A {\it conformal superalgebra} $R=R_0\oplus R_1$ is a $\Bbb{Z}_2$-graded $\Bbb{C}[\ptl]$-module with a $\Bbb{Z}_2$-graded linear map $Y^+(\cdot,z):\;R\rightarrow LM(R,R[z^{-1}]z^{-1})$ satisfying:
$$Y^+(\ptl u,z)={dY^+(u,z)\over dz}\qquad\for\;\;u\in R;\eqno(1.2)$$
$$Y^+(u,z)v=(-1)^{ij}\rd_x{e^{x\ptl}Y^+(v,-x)u\over z-x},\eqno(1.3)$$
$$Y^+(u,z_1)Y^+(v,z_2)-(-1)^{ij}Y^+(v,z_2)Y^+(u,z_1)=\rd_x{Y^+(Y^+(u,z_1-x)v,x)\over z_2-x} \eqno(1.4)$$
for $u\in R_i;\;v\in R_j$. We denote by $(R,\ptl,Y^+(\cdot,z))$ a conformal superalgebra. When $R_1=\{0\}$, we simply call $R$ a {\it conformal algebra}. 

The above definition is the equivalent generating-function form to that given in [Ka2], where the author used the component formulae with $Y^+(u,z)=\sum_{n=0}^{\infty}u_{(n)}z^{-1}$.

 Suppose that $(R,\ptl,Y^+(\cdot,z))$ is a conformal superalgebra that is a free $\Bbb{C}[\ptl]$-module over a $\Bbb{Z}_2$-graded subspace $V$, namely
$$R=\Bbb{C}[\ptl]V\;\;(\cong \Bbb{C}[\ptl]\otimes_{\Bbb{C}}V).\eqno(1.5)$$
Let $m$ be a positive integer. The algebra $R$ is called a {\it homogeneous conformal superalgebra of degree} $m$ if for any $u,v\in V$,
$$Y^+(u,z)v=\sum_{j=1}^m\ptl^{m-j}w_jz^{-j}\qquad\mbox{with}\;\;w_j\in V.\eqno(1.6)$$
A {\it  Lie superalgebra} $L$ is a $\Bbb{Z}_2$-graded algebra $L=L_0\oplus L_1$ with the operation $[\cdot,\cdot]$ satisfying
$$[u,v]=-(-1)^{ij}[v,u],\qquad [w,[v,u]]=[[w,v],u]-(-1)^{ij}[[w,u],v]\eqno(1.7)$$
for $u\in L_i,\;v\in L_j$ and $w\in L$. It is well known that a homogeneous conformal superalgebra of degree 1 is equivalent to a Lie superalgebra (e.g., cf. [Gdo1], [Ka2]). 

 A {\it Novikov superalgebra} is a ${\bf Z}_2$-graded vector space ${\cal A}={\cal A}_0\oplus {\cal A}_1$ with an operation ``$\circ$'' satisfying: for $u\in {\cal A}_i,\;v\in{\cal A}_j$ and $w\in {\cal A}_l$,
$$(u\circ v)\circ w=(-1)^{jl}(u\circ w)\circ v,\qquad (u,v,w)=(-1)^{ij}(v,u,w),\eqno(1.8)$$
where the associator:
$$(u,v,w)=(u\circ v)\circ w-u\circ (v\circ w).\eqno(1.9)$$
When ${\cal A}_1=\{0\}$, we call ${\cal A}$ a {\it Novikov algebra}. It was essentially stated in [GDo1] that a quadratic homogeneous conformal algebra is equivalent to a Novikov algebra. Such an algebraic structure appeared in [BN] from the point of view of Poisson structures of hydrodynamic type. The name ``Novikov algebra'' was given by Osborn [O1]. The above superanalogue was given in [X6]. 

By {\it quadratic conformal superalgebra}, we mean a conformal superalgebra $R$ that is a free $\Bbb{C}[\ptl]$-module over its $\Bbb{Z}_2$-graded subspace $V$ such that for $u,v\in V$,
$$Y^+(u,z)v=(w_1+\ptl w_2)z^{-1}+w_3z^{-2}\qquad \mbox{with}\;\;w_i\in V.\eqno(1.10)$$
It was essentially stated in [GDo1] (without proof) that a quadratic conformal superalgebra is equivalent to a bialgebraic structure $({\cal A},[\cdot,\cdot],\circ)$ such that $({\cal A},[\cdot,\cdot])$ forms a Lie algebra, $({\cal A},\circ)$ forms a Novikov algebra and the following compatibility condition holds:
$$[w\circ u,v]-[w\circ v,u]+[w,u]\circ v-[w,v]\circ u-w\circ [u,v]=0\eqno(1.11)$$
for $u,v,w\in {\cal A}$. We may call such a bialgebraic structure a {\it Gel'fand-Dorfman} bialgebra for convenience. A quadratic conformal algebra corresponds to a Hamiltonian pair in [GDo1], which plays fundamental roles in completely integrable systems. It was also pointed out in [GDo1] that if we define the commutator
$$[u,v]^-=u\circ v-v\circ u\eqno(1.12)$$
for a Novikov algebra $({\cal A},\circ)$, then $({\cal A},[\cdot,\cdot]^-,\circ)$ forms a Gel'fand-Dorfman bialgebra. 

In this paper, we shall study quadratic conformal superalgebras. Naturally, we need the following concept. A {\it Super Gel'fand-Dorfman bialgebra} is a $\Bbb{Z}_2$-graded vector space ${\cal A}={\cal A}_0\oplus {\cal A}_1$ with two algebraic operations $[\cdot,\cdot]$ and $\circ$ such that $({\cal A},[\cdot,\cdot])$ forms a Lie superalgebra, $({\cal A},\circ)$ forms a Novikov superalgebra and the following compatibility condition holds:
$$[w\circ u,v]-(-1)^{ij}[w\circ v,u]+[w,u]\circ v-(-1)^{ij}[w,v]\circ u-w\circ [u,v]=0\eqno(1.13)$$
for $u\in {\cal A}_i,\;v\in {\cal A}_j$ and $w\in {\cal A}$. We shall present the proof of that a quadratic conformal superalgebra is equivalent to a super Gel'fand-Dorfman bialgebra. For a Novikov superalgebra
$({\cal A},\circ)$, we define another operation $[\cdot,\cdot]^-$ on ${\cal A}$ by
$$[u,v]^-=u\circ v-(-1)^{ij}v\circ u\qquad\;\for\;\;u\in {\cal A}_i,\;v\in {\cal A}_j.\eqno(1.14)$$
The proof of the fact that $({\cal A},[\cdot,\cdot]^-,\circ)$ forms a super Gel'fand-Dorfman bialgebra will also be given. However, we do not claim these proofs as our major results in this paper. Our purpose of doing these is to give the reader a convenience.

Our main results can be divided into two aspects. First, we shall present four general constructions of super Gel'fand-Dorfman bialgebras. Two of them are extracted from simple Lie superalgebras of Cartan types W, H and K. One construction comes from a family of infinite-dimensional simple Lie superalgebras that we constructed in [X4]. The other construction is obtained from our classification works in this paper.

The second aspect of our results are classifications of Gel'fand-Dorfman bialgebras related to simple Novikov algebras. Zelmanov [Z] proved that any finite-dimensional simple Novikov algebra over an algebraically closed field with characteristic $0$ is one-dimensional. Osborn [O2] classified finite-dimensional simple Novikov algebras with an idempotent element over an algebraically closed field with prime characteristic. In [X3], we gave a complete classification of finite-dimensional simple Novikov algebras and their irreducible modules over an algebraically closed field with prime characteristic.   

Suppose that $\Delta$ is an additive subgroup of $\Bbb{C}$ and denote $\G=\{0\}$ or $\Bbb{N}=\{0,1,2,3,...\}$, the set of natural numbers. Let ${\cal A}_{\Delta,\G}$ be a vector space with a basis $\{x_{\al,j}\mid (\al,j)\in \Delta\times \G\}$. For any given constant $b\in \Bbb{C}$, we define algebraic operation $\circ$ on ${\cal A}_{\Delta,\G}$ by
$$x_{\al,i}\circ x_{\be,j}=(\be+b)x_{\al+\be,i+j}+jx_{\al+\be,i+j-1}\qquad\for\;\;\al,\be\in \Delta,\;i,j\in \G,\eqno(1.15)$$
where we adopt the convention that if a notion is not defined but technically appears in an expression, we always treat it as zero; for instance, $x_{\al,j}=0$ if $(\al,j)\not\in \Delta\times \G$ (this convention will be used throughout this paper). In [O3], Osborn proved that a simple Novikov algebras with an idempotent element whose left multiplication operator is locally-finite  over an algebraically closed field with  characteristic 0 must be isomorphic to $({\cal A}_{\Delta,\G},\circ)$ for some $\Delta,\;\G$ and $b$. There is a natural commutative associative algebra structure $\cdot$ on ${\cal A}_{\Delta,\G}$:
$$x_{\al,i}\cdot x_{\be,j}=x_{\al+\be,i+j}\qquad\for\;\;\al,\be\in \Delta,\;i,j\in\G.\eqno(1.16)$$
Throughout this paper, the symbol $\cdot$ of an associative algebraic operation in a product will be invisible for convenience, when the context is clear. For any $\xi\in {\cal A}_{\Delta,\G}$, we define an algebraic operation $\dmd_{\xi}$ on ${\cal A}_{\Delta,\G}$ by 
$$x_{\al,i}\dmd_{\xi} x_{\be,j}=(\be+\xi)x_{\al+\be,i+j}+jx_{\al+\be,i+j-1}\qquad\for\;\;\al,\be\in \Delta,\;i,j\in\G.\eqno(1.17)$$
We proved in [X4] that $({\cal A}_{\Delta,\G},\dmd_{\xi})$ forms a simple Novikov algebra for any $\xi\in {\cal A}_{\Delta,\G}$. This in particular gives a large family of simple Novikov algebras without any idempotent elements. According to Gel'fand and Dorfman's statement in [GDo1], we have a large family of Gel'fand-Dorfman bialgebra $({\cal A}_{\Delta,\G},[\cdot,\cdot]^-,\dmd_{\xi})$ (cf. (1.12)). In fact,
$$[x_{\al,i},x_{\be,j}]^-=x_{\al,i}\dmd_{\xi}x_{\be,j}-x_{\be,j}\dmd x_{\al,i}=(\be-\al)x_{\al+\be,i+j}+(j-i)x_{\al+\be,i+j-1}\eqno(1.18)$$
for $\al,\be\in \Delta,\;i,j\in\G$. Moreover, $({\cal A}_{\Delta,\G},[\cdot,\cdot]^-)$ is a simple Lie algebra (cf. [O5]). 

For convenience, we call $({\cal A},[\cdot,\cdot])$ a {\it Lie superalgebra over the Novikov superalgebra} $({\cal A},\circ)$ and $({\cal A},\circ)$ a {\it Novikov superalgebra over the Lie superalgebra} $({\cal A},[\cdot,\cdot])$ when $({\cal A},[\cdot,\cdot],\circ)$ forms a super Gel'fand-Dorfman bialgebra.
In [OZ], the authors proved that a Lie algebra over $({\cal A}_{\Delta,\G},\circ)$ with $\Delta=\Bbb{Z},\:\G=\{0\}$ and $b\not\in\Bbb{Z}$ or $\Delta=\{0\}$ and $\G=\Bbb{N}$ must be isomorphic to $({\cal A}_{\Delta,\G},[\cdot,\cdot]^-)$. In this paper, we shall classify all the Lie algebras over $({\cal A}_{\Delta,\{0\}},\circ)$ for arbitrary additive subgroup $\Delta$ and any constant $b$. In the case $b\not\in \Delta$, the Lie algebras are Block algebras [B]. We also classify the Lie algebras over $({\cal A}_{\Delta,\Bbb{N}},\circ)$ when $b\not\in\Delta$. It seems to us that there are too many complicated Lie algebras over $({\cal A}_{\Delta,\Bbb{N}},\circ)$ when $\Delta\neq \{0\}$ and $b\in\Delta$. We shall present several families of such Lie algebras which still look neat. Furthermore, we shall classify all the Novikov algebras whose commutator algebra is $({\cal A}_{\Delta,\{0\}},[\cdot,\cdot]^-)$. It looks more challenging to classify  Novikov superalgebras over  all the well-known simple Lie superalgebras. 

We remark that the Hamiltonian superoperator corresponding to a conformal superalgebra $R$ with $R_1\neq\{0\}$ does not in general have any analytic implication yet. The right theory of Hamiltonian operators compatible with supersymmetric partial differential equations is that we gave in [X6]. 

The paper is organized as follows. In Section 2, we mainly present the proof of the equivalence of a quadratic conformal superalgebra and a super Gel'fand-Dorfman  bialgebra. In Section 3, we give four general constructions of super Gel'fand-Dorfman bialgebras. Our classification results are presented in Sections 4, 5 and 6.

The author thanks Prof. E. Zel'manov for his comments and suggestions.

\section{Equivalence}

In this section, we shall prove that quadratic conformal superalgebras are equivalent to super Gel'fand-Dorfman bialgebras. Moreover, we shall also give the proof of that $({\cal A},[\cdot,\cdot]^-,\circ)$ forms a super Gel'fand-Dorfman bialgebra for any Novikov superalgebra $({\cal A},\circ)$ (cf. (1.14)).

Let $R$ be a $\Bbb{Z}_2$-graded free $\Bbb{C}[\ptl]$-module over its $\Bbb{Z}_2$-graded subspace $V$ $(\ptl(R_i)\subset R_i)$. Let $Y^+(\cdot,z):V\rightarrow LM(V,R[z^{-1}]z^{-1})$ be any given $\Bbb{Z}_2$-graded linear map. We can extend $Y^+(\cdot,z)$ as follows. First we extend $Y^+(\cdot,z)$ to a linear map from $R$ to $LM(V,R[z^{-1}])$:
$$Y^+(f(\ptl)u,z)v=f(d/dz)Y^+(u,z)v\qquad\for\;\;u,v\in V,\;f(\ptl)\in \Bbb{C}[\ptl].\eqno(2.1)$$
Then we define a linear map $Y^+(\cdot,z):R\rightarrow LM(R,R[z^{-1}])$ by
$$Y^+(\zeta,z)\ptl^mv=\sum_{j=0}^m(-1)^j\ptl^{m-j}d^jY^+(\zeta,z)v/dz^j\qquad\for\;\;\zeta\in R,\;v\in V,\;m\in \Bbb{N}.\eqno(2.2)$$
The extended map $Y^+(\cdot,z)$ naturally satisfies (1.2). According to Remark 4.1.2 (2) and (4.2.7) in [X7], $(R,Y^+(\cdot,z))$ forms a conformal superalgebra if and only if the map $Y^+(\cdot,z)$ satisfies (1.3) and (1.4) when acting on $V$ for $u,v\in V$. This fact was showed by Kac [Ka2] through a relatively more complicated approach. 

The map is said to be homogeneous of degree $m$ if it satisfies (1.6). Suppose that the map $Y^+(\cdot, z)$ can be written as:
$$Y^+(\cdot,z)=Y^+_1(\cdot,z)+Y_2^+(\cdot,z)+\cdots +Y^+_p(\cdot,z),\eqno(2.3)$$
where all $Y_j^+(\cdot,z):V\rightarrow LM(V,R[z^{-1}])$ are homogeneous $\Bbb{Z}_2$-graded linear maps of different degrees and they are supposed to be extended as $Y^+(\cdot,z)$. The following lemma can be proved by comparing the degrees in (1.3) and (1.4).
\psp

{\bf Lemma 2.1}. {\it The map} $Y^+(\cdot,z)$ {\it satisfies (1.3) if and only if all} $Y_j^+(\cdot,z)$ {\it satisfy (1.3). When} $p=2$, {\it the map} $Y^+(\cdot,z)$ {\it satisfies (1.4) if and only if} $Y_1(\cdot, z)$ {\it and} $Y_2(\cdot, z)$ {\it satisfy (1.4), and the following condition holds:}
\begin{eqnarray*}& &(Y^+_1(u,z_1)Y^+_2(v,z_2)+Y^+_2(u,z_1)Y^+_1(v,z_2))w\\& &
-(-1)^{ij}(Y^+_1(v,z_2)Y^+_2(u,z_1)+Y^+_2(v,z_2)Y^+_1(u,z_1))w\\&=&\rd_x{(Y_1^+(Y^+_2(u,z_1-x)v,x)+Y_2^+(Y^+_1(u,z_1-x)v,x))w\over z_2-x}\hspace{4cm}(2.4)\end{eqnarray*}
{\it for} $u\in V_i,\;v\in V_j$ {\it and} $w\in V$. {\it In particular, if the family} $(R,\ptl,Y^+(\cdot,z))$ {\it forms a conformal superalgebra, then both} $(R,\ptl,Y^+_1(\cdot,z))$ {\it and} $(R,\ptl,Y^+_2(\cdot,z))$ {\it form conformal superalgebras when} $p=2$.
\psp

{\bf Theorem 2.2}. {\it A quadratic conformal superalgebra is equivalent to a super Gel'fand-Dorfman bialgebra}.
\psp

{\it Proof}. Let $(R,\ptl,Y^+(\cdot,z))$ be a quadratic conformal superalgebra. By (1.10), we can write:
$$Y^+(u,z)v=(\ptl(v\circ u)+[v,u])z^{-1}+v\cdot uz^{-2}\qquad\for\;\;u,v\in V,\eqno(2.5)$$
where $\circ,\;[\cdot,\cdot],\;\cdot$ are three algebraic operations on $V$ and the reason of changing the order of $u$ and $v$ on the right hand side is because we want our results consistent with the notions of Novikov algebra and Gel'fand-Dorfman bialgebra (otherwise, we would obtain their ``opposite algebras'').
For $u\in V_i$ and $v\in V_j$, (1.3) becomes:
\begin{eqnarray*}& &(\ptl(u\circ v)+[u,v])z^{-1}+u\cdot vz^{-2} \\&=&(-1)^{ij}[(\ptl(v\cdot u-v\circ u)-[v,u])z^{-1}+v\cdot uz^{-2}],\hspace{5.7cm}(2.6)\end{eqnarray*}
where $f(x)\in R[x]$. Since $R$ is a free $\Bbb{C}[\ptl]$-module over $V$, we obtain:
$$v\cdot u=v\circ u+(-1)^{ij}u\circ v,\qquad [u,v]=-(-1)^{ij}[v,u].\eqno(2.7)$$

Define two linear maps $Y_1^+(\cdot,z),\;Y_2^+(\cdot,z):V\rightarrow LM(V,R[z^{-1}])$ by:
$$Y_1^+(u,z)v=[v,u]z^{-1},\;\;Y_2^+(u,z)v=\ptl v\circ uz^{-1}+(v\circ u+(-1)^{ij}u\circ v)z^{-2}\eqno(2.8)$$
for $u\in V_i,\;v\in V_j$. Then $Y_j^+(\cdot,z)$ is a homogeneous map of degree $j$ and  $Y^+(\cdot,z)=Y_1^+(\cdot,z)+Y_2^+(\cdot,z)$. 
According to the above lemma, both $Y_1^+(\cdot,z)$ and $Y_2^+(\cdot,z)$ satisfy (1.4). Let $u\in V_i,\;v\in V_j$. By (1.4) for $Y_1^+(\cdot,z)$, we have:
$$([[w,v],u]-(-1)^{ij}[[w,u],v])z_1^{-1}z_2^{-1}=[w,[v,u]]z^{-1}_1z_2^{-1}\eqno(2.9)$$
by Definition 2.7b in [Ka2], which implies the second equation in (1.7). Thus $(V,[\cdot,\cdot])$ forms a Lie superalgebra. 

Next we let $u\in V_i,\;v\in V_j$ and $w\in V_l$. By (2.2) and Definition 2.7b in [Ka2],
\begin{eqnarray*}& &Y_2^+(u,z_1)Y_2^+(v,z_2)w\\&=&\ptl^2(w\circ v)\circ uz_1^{-1}z_2^{-1}+\ptl((w\circ v)\circ u+(-1)^{jl}(v\circ w)\circ u)z_1^{-1}z_2^{-2}\\& &+\ptl[2(w\circ v)\circ u+(-1)^{i(j+l)}u\circ (w\circ v)]z_1^{-2}z_2^{-1}\\&&+[(w\circ v+(-1)^{jl}v\circ w)\circ u+(-1)^{i(j+l)}u\circ (w\circ v+(-1)^{jl}v\circ w)]z_1^{-2}z_2^{-2}\\& &+2((w\circ v)\circ u+(-1)^{i(j+l)}u\circ (w\circ v))z_1^{-3}z_2^{-1},\hspace{5.8cm}(2.10)\end{eqnarray*}
\begin{eqnarray*}& &\rd_x{Y^+_2(Y^+_2(u,z_1-x)v,x)w\over z_2-x}\\&=&-\ptl w\circ (v\circ u)z_1^{-1}z_2^{-2}-2(w\circ(v\circ u)+(-1)^{l(i+j)}(v\circ u)\circ w)z_1^{-1}z_2^{-3}+(-w\circ(v\circ u)\\& &+(-1)^{ij}w\circ(u\circ v)-(-1)^{l(i+j)}(v\circ u)\circ w+(-1)^{l(i+j)+ij}(u\circ v)\circ w)z_1^{-2}z_2^{-2}\\& &+(-1)^{ij}\ptl w\circ(u\circ v)z_1^{-2}z_2^{-1}+(-1)^{ij}2(w\circ(u\circ v)\\& &+(-1)^{l(i+j)}(u\circ v)\circ w)z_1^{-3}z_2^{-1}\hspace{8.6cm}(2.11)\end{eqnarray*}
by (2.1). Thus by (1.4),
\begin{eqnarray*}& &\ptl^2[(w\circ v)\circ u-(-1)^{ij}(w\circ u)\circ v]z_1^{-1}z_2^{-1}-(-1)^{ij}2((w\circ u)\circ v\\& &+(-1)^{j(i+l)}v\circ (w\circ u))z_1^{-1}z_2^{-3}+\ptl[(w\circ v)\circ u+(-1)^{jl}(v\circ w)\circ u\\& &-(-1)^{ij}(2(w\circ u)\circ v+(-1)^{j(i+l)}v\circ (w\circ u))]z_1^{-1}z_2^{-2}-(-1)^{ij}\ptl[(w\circ u)\circ v\hspace{4cm}\end{eqnarray*}
\begin{eqnarray*}& &+(-1)^{il}(u\circ w)\circ v-(-1)^{ij}(2(w\circ v)\circ u+(-1)^{i(j+l)}u\circ (w\circ v))]z_1^{-2}z_2^{-1}\\& &+\{[(w\circ v+(-1)^{jl}v\circ w)\circ u+(-1)^{i(j+l)}u\circ (w\circ v+(-1)^{jl}v\circ w)]\\& &-(-1)^{ij}[(w\circ u+(-1)^{il}u\circ w)\circ v+(-1)^{j(i+l)}v\circ (w\circ u\\& &+(-1)^{il}u\circ w)]\}z_1^{-2}z_2^{-2}+2((w\circ v)\circ u+(-1)^{i(j+l)}u\circ (w\circ v))z_1^{-3}z_2^{-1}\\&=&Y_2^+(u,z_1)Y_2^+(v,z_2)w-(-1)^{ij}Y_2^+(v,z_2)Y_2^+(u,z_1)w\\&=&\rd_x{Y^+_2(Y^+_2(u,z_1-x)v,x)w\over z_2-x}\\&=&-\ptl w\circ (v\circ u)z_1^{-1}z_2^{-2}-2(w\circ(v\circ u)+(-1)^{l(i+j)}(v\circ u)\circ w)z_1^{-1}z_2^{-3}+(-w\circ(v\circ u)\\& &+(-1)^{ij}w\circ(u\circ v)-(-1)^{l(i+j)}(v\circ u)\circ w+(-1)^{l(i+j)+ij}(u\circ v)\circ w)z_1^{-2}z_2^{-2}\\& &+(-1)^{ij}\ptl w\circ(u\circ v)z_1^{-2}z_2^{-1}+(-1)^{ij}2(w\circ(u\circ v)\\& &+(-1)^{l(i+j)}(u\circ v)\circ w)z_1^{-3}z_2^{-1}.\hspace{8.6cm}(2.12)\end{eqnarray*}

Comparing the coefficients of $z_1^{-1}z_2^{-1}$ in (2.12), we have:
$$(w\circ v)\circ u=(-1)^{ij}(w\circ u)\circ v.\eqno(2.13)$$
From the coefficients of $z_1^{-1}z_2^{-2}$ in (2.12), we find
\begin{eqnarray*}& &(w\circ v)\circ u+(-1)^{jl}(v\circ w)\circ u-(-1)^{ij}(2(w\circ u)\circ v+(-1)^{j(i+l)}v\circ (w\circ u))\\&=&- w\circ (v\circ u),\hspace{11.5cm}(2.14)\end{eqnarray*}
which is equivalent to:
$$(-1)^{jl}(v\circ w)\circ u-(w\circ v)\circ u-(-1)^{jl}v\circ (w\circ u)=-w\circ (v\circ u)\eqno(2.15)$$
by (2.13). Note that (2.15) can be written as:
$$(w,v,u)=(-1)^{jl}(v,w,u)\eqno(2.16)$$
(cf. (1.9)). Hence $(V,\circ)$ forms a Novikov superalgebra. 

We want to show that the coefficients of the other monomials in (2.12) do not give more constraints on $(V,\circ)$. 
Checking the coefficient of $z_1^{-1}z_2^{-3}$, we have:
$$-(-1)^{ij}2((w\circ u)\circ v+(-1)^{j(i+l)}v\circ (w\circ u))=-2(w\circ(v\circ u)+(-1)^{l(i+j)}(v\circ u)\circ w),\eqno(2.17)$$
which is equivalent to:
$$(w\circ v)\circ u+(-1)^{jl}v\circ (w\circ u)=w\circ(v\circ u)+(-1)^{lj}(v\circ w)\circ u\eqno(2.18)$$
by (2.13). Note (2.18) is the same as (2.15).
The equation from the coefficient of $z_1^{-2}z_2^{-1}$ is equivalent to that from the coefficient of $z_1^{-1}z_2^{-2}$ and
the equation from the coefficient of $z_1^{-3}z_2^{-1}$ is equivalent to that from the coefficient of $z_1^{-1}z_2^{-3}$ because $u,v,w$ are arbitrary.
Extracting the coefficient of $z_1^{-2}z_2^{-2}$ in (2.12), we obtain:
\newpage

\begin{eqnarray*}& &(w\circ v+(-1)^{jl}v\circ w)\circ u+(-1)^{i(j+l)}u\circ (w\circ v+(-1)^{jl}v\circ w)\\& &-(-1)^{ij}[(w\circ u+(-1)^{il}u\circ w)\circ v+(-1)^{j(i+l)}v\circ (w\circ u+(-1)^{il}u\circ w)]\\&=&-w\circ(v\circ u)+(-1)^{ij}w\circ(u\circ v)-(-1)^{l(i+j)}(v\circ u)\circ w\\& &+(-1)^{l(i+j)+ij}(u\circ v)\circ w,\hspace{9.4cm}(2.19)\end{eqnarray*}
which is equivalent to:
\begin{eqnarray*}& &(-1)^{jl}(v,w,u)-(-1)^{i(j+l)}(u,w,v)\\&=&(w,v,u)-(-1)^{ij}(w,u,v)-(-1)^{l(i+j)}(v,u,w)+(-1)^{ij+l(i+j)}(u,v, w).\hspace{1.4cm}(2.20)\end{eqnarray*}
The above equation is implied by (2.16) again because $u,v,w$ are arbitrary.

Now we want to show that (2.4) only give rise to (1.13). First we have:
\begin{eqnarray*}& &(Y_1^+(u,z_1)Y_2^+(v,z_2)+Y_2^+(u,z_1)Y_1^+(v,z_2))w\\&=&\ptl ([w\circ v,u]+[w,v]\circ u)z_1^{-1}z_2^{-1}+[w\circ v+(-1)^{jl}v\circ w,u]z_1^{-1}z_2^{-2}\\&&+([w,v]\circ u+(-1)^{i(j+l)}u\circ [w,v]+[w\circ v,u])z_1^{-2}z_2^{-1},\hspace{4.5cm}(2.21)\end{eqnarray*}
\begin{eqnarray*}&&\rd_x(z_2-x)^{-1}[Y_1^+(Y^+_2(u,z_1-x)v,x)w+Y_2^+(Y^+_1(u,z_1-x)v,x)w]\\&=&\ptl w\circ [v,u]z_1^{-1}z_2^{-1}+(w\circ[v,u]+(-1)^{l(i+j)}[v,u]\circ w-[w,v\circ u])z_1^{-1}z_2^{-2}\\& &+(w\circ[v,u]+(-1)^{l(i+j)}[v,u]\circ w+(-1)^{ij}[w,u\circ v])z_1^{-2}z_2^{-1}.\hspace{3.5cm}(2.22)\end{eqnarray*}
Then by (2.4), we get
\begin{eqnarray*}& &\ptl ([w\circ v,u]+[w,v]\circ u-(-1)^{ij}([w\circ u,v]+[w,u]\circ v))z_1^{-1}z_2^{-1}\\& &
+([w\circ v+(-1)^{jl}v\circ w,u]-(-1)^{ij}([w,u]\circ v+(-1)^{j(i+l)}v\circ [w,u]+[w\circ u,v]))z_1^{-1}z_2^{-2}\\&&+([w,v]\circ u+(-1)^{i(j+l)}u\circ [w,v]+[w\circ v,u]-(-1)^{ij}[w\circ u+(-1)^{jl}u\circ w,v])z_1^{-2}z_2^{-1}\\&=&(Y_1^+(u,z_1)Y_2^+(v,z_2)+Y_2^+(u,z_1)Y_1^+(v,z_2))w\\&&-(-1)^{ij}(Y_1^+(v,z_2)Y_2^+(u,z_1)+Y_2^+(v,z_1)Y_1^+(u,z_1))w\\&=&\rd_x(z_2-x)^{-1}[Y_1^+(Y^+_2(u,z_1-x)v,x)w+Y_2^+(Y^+_1(u,z_1-x)v,x)w]\\&=&\ptl w\circ [v,u]z_1^{-1}z_2^{-1}+(w\circ[v,u]+(-1)^{l(i+j)}[v,u]\circ w-[w,v\circ u])z_1^{-1}z_2^{-2}\\& &+(w\circ[v,u]+(-1)^{l(i+j)}[v,u]\circ w+(-1)^{ij}[w,u\circ v])z_1^{-2}z_2^{-1}.\hspace{3.4cm}(2.23)\end{eqnarray*}

Note that the coefficients of $z_1^{-1}z_2^{-1}$ in (2.23) imply:
$$[w\circ v,u]+[w,v]\circ u-(-1)^{ij}([w\circ u,v]+[w,u]\circ v)=w\circ [v,u],\eqno(2.24)$$
which is the same as (1.13) because $u,v,w$ are arbitrary.
Comparing the coefficients of $z_1^{-1}z_2^{-2}$ in (2.23), we obtain:
\newpage

\begin{eqnarray*}& &[w\circ v+(-1)^{jl}v\circ w,u]-(-1)^{ij}([w,u]\circ v+(-1)^{j(i+l)}v\circ [w,u]+[w\circ u,v])\\&=&w\circ[v,u]+(-1)^{l(i+j)}[v,u]\circ w-[w,v\circ u],\hspace{6.3cm}(2.25)\end{eqnarray*}
which can be rewritten as:
\begin{eqnarray*}& &([w\circ v,u]-(-1)^{ij}[w\circ u,v]-(-1)^{ij}[w,u]\circ v-w\circ[v,u]+[w,v]\circ u)\\& &+((-1)^{jl}[v,w]\circ u+(-1)^{jl}[v\circ w,u]-(-1)^{l(i+j)}[v\circ u,w]\\& &-(-1)^{l(i+j)}[v,u]\circ w-(-1)^{jl}v\circ [w,u])=0.\hspace{6.4cm}(2.26)\end{eqnarray*}
However, (2.26) is implied by (2.24) because $u,v,w$ are arbitrary. The equation from the coefficients of $z_1^{-2}z_2^{-1}$ in (2.23) is equivalent to (2.25), again because $u,v,w$ are arbitrary. 

What we have proved in the above is that $(V,[\cdot,\cdot],\circ)$ forms a super Gel'fand-Dorman bialgebra. Since the above arguments are reversible, a super Gel'fand-Dorman bialgebra $(V,[\cdot,\cdot],\circ)$  defines
a quadratic conformal superalgebra by (2.5) and the first equation in (2.7). Thus the equivalence is established.$\qquad\Box$
\psp

{\bf Theorem 2.3}. {\it Let} $({\cal A},\circ)$ {\it be a Novikov superalgebra}. {\it Then} $({\cal A},[\cdot,\cdot]^-,\circ)$ {\it (cf. (1.14)) forms a super Gel'fand-Dorfman bialgebra}.
\psp

{\it Proof}. For $u\in {\cal A}_i,\;v\in {\cal A}_j$ and $w\in {\cal A}_l$,
\begin{eqnarray*}& &[w\circ u,v]^--(-1)^{ij}[w\circ v,u]^-+[w,u]^-\circ v-(-1)^{ij}[w,v]^-\circ u-w\circ[u,v]^-\\&=&(w\circ u)\circ v-(-1)^{j(i+l)}v\circ(w\circ u)-(-1)^{ij}(w\circ v)\circ u+(-1)^{il}u\circ(w\circ v)\\&&+(w\circ u)\circ v-(-1)^{il}(u\circ w)\circ v-(-1)^{ij}(w\circ v)\circ u\\& &+(-1)^{(i+l)j}(v\circ w)\circ u-w\circ(u\circ v)+(-1)^{ij}w\circ(v\circ u)\\&=&
-(-1)^{j(i+l)}v\circ(w\circ u)+(-1)^{il}u\circ(w\circ v)-(-1)^{il}(u\circ w)\circ v\\& &+(-1)^{(i+l)j}(v\circ w)\circ u-w\circ(u\circ v)+(-1)^{ij}w\circ(v\circ u)\\&=&(-1)^{(i+l)j}(v,w,u)-(-1)^{il}(u,w,v)-w\circ(u\circ v)+(-1)^{ij}w\circ(v\circ u)\\&=&(-1)^{ij}(w,v,u)-(w,u,v)-w\circ(u\circ v)+(-1)^{ij}w\circ(v\circ u)\\&=&
(-1)^{ij}(w\circ v)\circ u-(-1)^{ij}w\circ (v\circ u)-(w\circ u)\circ v+w\circ (u\circ v)\\&&-w\circ(u\circ v)+(-1)^{ij}w\circ(v\circ u)=0\hspace{7.3cm}(2.27)\end{eqnarray*}
by (1.8). So (1.13) holds.$\qquad\Box$
\psp

{\bf Example}. A {\it super commutative associative algebra} ${\cal A}$ is $\Bbb{Z}_2$-graded associative algebra such that
\newpage

$$u\cdot v=(-1)^{ij}v\cdot u\qquad\for\;\;u\in {\cal A}_i,\;v\in {\cal A}_j.\eqno(2.28)$$
It is easily verified that a super commutative associative algebra forms a Novikov superalgebra. 

An element $d\in \mbox{End}\:{\cal A}$ is called a {\it derivation} of ${\cal A}$ if there exists $i\in \Bbb{Z}_2$ such that
$$d({\cal A}_j)\subset {\cal A}_{i+j},\qquad d(u\cdot v)=d(u)\cdot v+(-1)^{ij}v\cdot d(v),\eqno(2.29)$$
for $j\in \Bbb{Z}_2,\;u\in {\cal A}_j,\;v\in{\cal A}.$ The derivation $d$ is called {\it even} if $i=0$ and is called {\it odd} if $i=1$. 

Let $d$ be an even derivation of a super commutative associative algebra ${\cal A}$ and $\xi\in {\cal A}_0$. We define an operation $\circ$ on ${\cal A}$ by
$$u\circ v=ud(v)+\xi u v\qquad\for\;\;u,v\in{\cal A}.\eqno(2.30)$$
By the analogous arguments as (2.7) and (2.8) in [X4], we can prove that $({\cal A},\circ)$ forms a Novikov superalgebra.

\section{Constructions}

In this section, we shall give four constructions of super Gel'fand-Dorfman bialgebras.

Let $({\cal A},\cdot)$ be a super commutative associative algebra. Denote by $W({\cal A})_0$ the space of even derivations of ${\cal A}$ and by $W({\cal A})_1$ the space of odd derivations of ${\cal A}$. Then
$$ W({\cal A})=W({\cal A})_0+W({\cal A})_1\eqno(3.1)$$
forms a Lie superalgebra with respect to $[\cdot,\cdot]$ defined by
$$[d_1,d_2](v)=d_1(d_2(v))-(-1)^{ij}d_2(d_1(v))\qquad\for\;\;d_1\in W({\cal A})_i,\;d_2\in W({\cal A})_j,\;v\in {\cal A}.\eqno(3.2)$$
Moreover, $W({\cal A})$ forms a left ${\cal A}$-module with respect to the action:
$$(ad)(v)=a d(v)\qquad\for\;\;a,v\in {\cal A},\;d\in W({\cal A}).\eqno(3.3)$$
Set
$${\cal N}=W({\cal A})\oplus {\cal A}.\eqno(3.4)$$
We define two algebraic operation $[\cdot,\cdot]$ and $\circ$ on ${\cal N}$ by
$$[d_1+\xi_1,d_2+\xi_2]=[d_1,d_2]+d_1(\xi_2)-(-1)^{ij}d_2(\xi_1),\eqno(3.5)$$
$$(d_1+\xi_1)\circ(d_2+\xi_2)=(-1)^{ij}\xi_2d_1+\xi_1\xi_2\eqno(3.6)$$
for $d_1+\xi_1\in {\cal N}_i=W({\cal A})_i+{\cal A}_i$ and $d_2+\xi_2\in {\cal N}_j$. 
\psp

{\bf Theorem 3.1}. {\it The family} $({\cal N},[\cdot,\cdot],\circ)$ {\it forms a super Gel'fand-Dorfman bialgebra}.
\psp

{\it Proof}. The pair $({\cal N},[\cdot,\cdot])$ forms a Lie superalgebra because it is a semi-product of $W({\cal A})$ with its module ${\cal A}$. Let $d_i+\xi_i\in {\cal N}_{j_i}$ with $i=1,2,3$. First we have:
\begin{eqnarray*}[(d_1+\xi_1)\circ (d_2+\xi_2)]\circ (d_3+\xi_3)&=&(-1)^{j_1j_2+(j_1+j_2)j_3}\xi_3\xi_2d_1+\xi_1\xi_2\xi_3\\&=&(-1)^{j_2j_3}[(d_1+\xi_1)\circ (d_3+\xi_3)]\circ (d_2+\xi_2),\hspace{1.1cm}(3.7)\end{eqnarray*}
\begin{eqnarray*}(d_1+\xi_1)\circ [(d_2+\xi_2)\circ (d_3+\xi_3)]&=&(d_1+\xi_1)\circ ((-1)^{j_2j_3}\xi_3d_2+\xi_2\xi_3)\\&=&(-1)^{j_1(j_2+j_3)}\xi_2\xi_3d_1+\xi_1\xi_2\xi_3.\hspace{3.5cm}(3.8)\end{eqnarray*}
The above two expressions imply the associativity:
$$[(d_1+\xi_1)\circ (d_2+\xi_2)]\circ (d_3+\xi_3)=(d_1+\xi_1)\circ [(d_2+\xi_2)\circ (d_3+\xi_3)].\eqno(3.9)$$
Thus $({\cal N},\circ)$ forms a Novikov superalgebra. Furthermore,
\begin{eqnarray*}& &[(d_3+\xi_3)\circ (d_1+\xi_1),d_2+\xi_2]-(-1)^{j_1j_2}[(d_3+\xi_3)\circ (d_2+\xi_2),d_1+\xi_1]\\& &+[d_3+\xi_3,d_1+\xi_1]\circ (d_2+\xi_2)-(-1)^{j_1j_2}[d_3+\xi_3,d_2+\xi_2]\circ (d_1+\xi_1)\hspace{3cm}\end{eqnarray*}
\begin{eqnarray*}& &-(d_3+\xi_3)\circ [d_1+\xi_1,d_2+\xi_2]\\&=&[(-1)^{j_1j_3}\xi_1d_3+\xi_3\xi_1,d_2+\xi_2]-(-1)^{j_1j_2}[(-1)^{j_2j_3}\xi_2d_3+\xi_3\xi_2,d_1+\xi_1]\\& &+([d_3,d_1]+d_3(\xi_1)-(-1)^{j_1j_3}d_1(\xi_3))\circ (d_2+\xi_2)-(-1)^{j_1j_2}([d_3,d_2]+d_3(\xi_2)\\& &-(-1)^{j_2j_3}d_2(\xi_3))\circ (d_1+\xi_1)-(d_3+\xi_3)\circ ([d_1,d_2]+d_1(\xi_2)-(-1)^{j_1j_2}d_2(\xi_1))\\&=&(-1)^{j_1j_3}\xi_1[d_3,d_2]-(-1)^{j_1j_3+j_2(j_1+j_3)}d_2(\xi_1)d_3+(-1)^{j_1j_3}\xi_1d_3(\xi_2)\\& &-(-1)^{j_2(j_1+j_3)}d_2(\xi_3\xi_1)-(-1)^{j_2(j_1+j_3)}\xi_2[d_3,d_1]+(-1)^{(j_1+j_2)j_3}d_1(\xi_2)d_3\\& &-(-1)^{j_2(j_1+j_3)}\xi_2d_3(\xi_1)+(-1)^{j_1j_3}d_1(\xi_3\xi_2)+(-1)^{j_2(j_1+j_3)}\xi_2[d_3,d_1]\\& &+(d_3(\xi_1)-(-1)^{j_1j_3}d_1(\xi_3))\xi_2-(-1)^{j_1j_3}\xi_1[d_3,d_2]-(-1)^{j_1j_2}(d_3(\xi_2)\\& &-(-1)^{j_2j_3}d_2(\xi_3))\xi_1-(-1)^{j_3(j_1+j_2)} d_1(\xi_2)d_3+(-1)^{j_1j_2+j_3(j_1+j_2)}d_2(\xi_1)d_3\\&&-\xi_3d_1(\xi_2)+(-1)^{j_1j_2}\xi_3d_2(\xi_1)=0.\qquad\Box\hspace{6.8cm}(3.10)\end{eqnarray*}

The above construction is extracted from the simple Lie superalgebras of Cartan type $W$ (cf. [Ka3], [X7]). 

Our second and third constructions are related to the following concept. A {\it Lie-Poisson superalgebra} ${\cal A}$ is a $\Bbb{Z}_2$-graded space with two algebraic operations $\cdot$ and $[\cdot,\cdot]$ such that $({\cal A},\cdot)$ forms a super commutative associative algebra, $({\cal A},[\cdot,\cdot])$ forms a Lie superalgebra and the following compatibility condition is satisfied:
$$[u,v\cdot w]=[u,v]\cdot w+(-1)^{ij}v\cdot [u,w]\qquad\for\;\;u\in{\cal A}_i,\;v\in{\cal A}_j,\;w\in{\cal A}.\eqno(3.11)$$

 Let $({\cal A},\cdot,[\cdot,\cdot])$ be a Lie-Poisson superalgebra and let $d$ be an even derivation of the algebra $({\cal A},\cdot)$ such that
$$d[u,v]=[d(u),v]+[u,d(v)]+\xi[u,v]\qquad\for\;\;u,v\in{\cal A},\eqno(3.12)$$
where $\xi\in \Bbb{F}$ is a constant.
Now we define another algebraic operation on $\circ$ on ${\cal A}$ by
$$u\circ v=ud(v)+\xi uv\qquad\for\;\;u\in {\cal A},\;v\in {\cal A}_{\be}.\eqno(3.13)$$

{\bf Theorem 3.2}. {\it The family} $({\cal A},[\cdot,\cdot],\circ)$ {\it forms a Gel'fand-Dorfman super bialgebra}.
\psp

{\it Proof}. Note that (3.13) is the same as the second equation in (2.15) in [X4]. It can be similarly verified that $({\cal A},\circ)$ forms a Novikov superalgebra. Moreover, for $u\in{\cal A}_i,\;v\in {\cal A}_j$ and $w\in {\cal A}$,
\begin{eqnarray*}& &[w\circ u,v]-(-1)^{ij}[w\circ v,u]+[w,u]\circ v-(-1)^{ij}[w,v]\circ u-w\circ [u,v]\\&=&
[w(d+\xi)(u),v]-(-1)^{ij}[w(d+\xi)(v),u]+[w,u](d+\xi)(v)\\& &-(-1)^{ij}[w,v](d+\xi)(u)-w(d+\xi) [u,v]\\&=&
[wd(u),v]-(-1)^{ij}[wd(v),u]+[w,u]d(v)-(-1)^{ij}[w,v]d(u)+\xi([wu,v]\\& &-(-1)^{ij}[wv,u]+[w,u]v-(-1)^{ij}[w,v]u)
-w([d(u),v]+[u,d(v)]+2\xi[u,v])\\&=&(-1)^{ij}[w,v]d(u)-[w,u]d(v)+[w,u]d(v)-(-1)^{ij}[w,v]d(u)\\& &+\xi([w[u,v]-(-1)^{ij}w[v,u])-w(2\xi[u,v])=0.\hspace{5.4cm}(3.14)\end{eqnarray*}
So (1.13) holds.$\qquad\Box$
\psp

The above construction is related to the Lie superalgebras of Hamiltonian type and Contact type (cf. [Ka1]).

Next we shall present a construction extracted from a family of infinite-dimensional simple Lie superalgebras that we obtained in [X4] (cf. Theorem 5.3 in [X4])).

Let $({\cal A},\cdot)$ be a commutative associative algebra and let $d$ be a derivation of ${\cal A}$. Set
$$\td{\cal A}={\cal A}\times {\cal A}=\td{\cal A}_0\oplus \td{\cal A}_1\eqno(3.15)$$
with
$$\td{\cal A}_0=({\cal A},0),\qquad \td{\cal A}_1=(0,{\cal A}).\eqno(3.16)$$
 For fixed elements $\xi,\eta_0,\eta_1\in {\cal A}$, we define two algebraic operations $[\cdot,\cdot]$ and $\circ$ on $\td{\cal A}$ by
$$[(u_0,u_1),(v_0,v_1)]=(\xi u_1v_1,0),\eqno(3.17)$$
$$(u_0,u_1)\circ(v_0,v_1)=(u_0(d(v_0)+\eta_0v_0),u_1(d(v_0)+\eta_0v_0)+u_0(d(v_1)+\eta_1v_1))\eqno(3.18)$$
for $u_i,v_j\in {\cal A}$. It is easily seen that $(\td{\cal A},[\cdot,\cdot])$ forms a Lie superalgebra.

First we want to show that $(\td{\cal A},\circ)$ forms a Novikov superalgebra. Note that
by Corollary 2.6 in [X4], $(\td{\cal A}_0,\circ)$ forms a Novikov algebra. This fact will simplify our verification of (1.13). Let $u_i,v_i,w_i\in {\cal A}$ with $i=0,1$. First we have:
\begin{eqnarray*}[(0,w_1)\circ (u_0,0)]\circ (v_0,0)&=&(0,w_1(d(u_0)+\eta_0u_0))\circ (v_0,0)\\&=&(0,w_1(d(u_0)+\eta_0u_0)(d(v_0)+\eta_0v_0)),\hspace{3.1cm}(3.19)\end{eqnarray*}
\begin{eqnarray*}[(w_0,w_1)\circ (0,u_1)]\circ (v_0,0)&=&(0,w_0(d(u_1)+\eta_1u_1))\circ (v_0,0)\\&=&(0,w_0(d(u_1)+\eta_1u_1)(d(v_0)+\eta_0v_0)),\hspace{2.8cm}(3.20)\end{eqnarray*}
\begin{eqnarray*}[(w_0,w_1)\circ (v_0,0)]\circ (0,u_1)&=&(w_0(d(v_0)+\eta_0v_0),w_1(d(v_0)+\eta_0v_0))\circ (0,u_1)\\&=&(0,w_0(d(v_0)+\eta_0v_0)(d(u_1)+\eta_1u_1)),\hspace{2.9cm}(3.21)\end{eqnarray*}
$$[(w_0,w_1)\circ (0,u_1)]\circ (0,v_1)=(0,w_0(d(u_1)+\eta_1u_1))\circ (0,v_1)=(0,0).\eqno(3.22)$$
Thus the first equation in (1.8) holds. Furthermore,
\begin{eqnarray*}& &((u_0,0),(v_0,0),(0,w_1))\\&=&[(u_0,0)\circ(v_0,0)]\circ(0,w_1)-(u_0,0)\circ[(v_0,0)\circ(0,w_1)]\\&=&(0,u_0(d(v_0)+\eta_0v_0)(d(w_1)+\eta_1w_1)-(0,u_0(d(v_0d(w_1)+\eta_1v_0w_1)\\& &+\eta_1(v_0d(w_1)+\eta_1v_0w_1)))\\&=&(0, u_0v_0[(\eta_0-\eta_1)(d(w_1)+\eta_1w_1)-d^2(w_1)-d(\eta_1w_1)]),\hspace{4.2cm}(3.23)\end{eqnarray*}
\begin{eqnarray*}& &((u_0,0),(0,v_1),(w_0,w_1))\\&=&[(u_0,0)\circ(0,v_1)]\circ(w_0,w_1)-(u_0,0)\circ[(0,v_1)\circ(w_0,w_1)]\\&=&(0,u_0(d(v_1)+\eta_1v_1)(d(w_0)+\eta_0w_0))-(0,u_0(d(v_1d(w_0)+\eta_0v_1w_0)\\& &+\eta_1(v_1d(w_0)+\eta_0v_1w_0)))\\&=&(0,-u_0v_1d(d(w_0)+\eta_0w_0)),\hspace{9cm}(3.24)\end{eqnarray*}
\begin{eqnarray*}& &((0,v_1),(u_0,0),(w_0,w_1))\\&=&[(0,v_1)\circ(u_0,0)]\circ(w_0,w_1)-(0,v_1)\circ[(u_0,0)\circ(w_0,w_1)]\\&=&(0,v_1(d(u_0)+\eta_0u_0)(d(w_0)+\eta_0w_0))-(0,v_1(d(u_0d(w_0)+\eta_0u_0w_0)\\&&+\eta_0(u_0d(w_0)+\eta_0u_0w_0)))\\&=&(0,-u_0v_1d(d(w_0)+\eta_0w_0)),\hspace{9cm}(3.25)\end{eqnarray*}
$$((0,u_1),(0,v_1),(w_0,w_1))=0=-((0,v_1),(0,u_1),(w_0,w_1)).\eqno(3.26)$$
Expressions (3.23-26) show that the second equation in (1.8) holds. Hence $(\td{\cal A},\circ)$ forms a Novikov superalgebra. 

We shall now exam (1.13). Note (1.13) trivially holds if $u=(u_0,0)$ and $v=(v_0,0)$ by (3.17). Moreover,
\begin{eqnarray*}& &[(w_0,w_1)\circ (u_0,0),(0,v_1)]-[(w_0,w_1)\circ (0,v_1),(u_0,0)]+[(w_0,w_1),(u_0,0)]\circ(0,v_1)\\& &-[(w_0,w_1),(0,v_1)]\circ(u_0,0)-(w_0,w_1)\circ[(u_0,0),(0,v_1)]\\&=&(\xi w_1v_1(d(u_0)+\eta_0u_0),0)-0+0-(\xi w_1v_1(d(u_0)+\eta_0u_0),0)-0=0,\hspace{1.8cm}(3.27)\end{eqnarray*}
\begin{eqnarray*}& &[(w_0,w_1)\circ (0,u_1),(0,v_1)]+[(w_0,w_1)\circ (0,v_1),(0,u_1)]\\&&+[(w_0,w_1),(0,u_1)]\circ(0,v_1)+[(w_0,w_1),(0,v_1)]\circ(0,u_1)-(w_0,w_1)\circ[(0,u_1),(0,v_1)]\\&=&(\xi w_0v_1(d(u_1)+\eta_1 u_1),0)+(\xi w_0u_1(d(v_1)+\eta_1 v_1),0)+(0,\xi w_1u_1(d(v_1)+\eta_1 v_1))\\&& +(0,\xi w_1v_1(d(u_1)+\eta_1 u_1))-(w_0(d(\xi u_1v_1)+\eta_0\xi u_1v_1),w_1(d(\xi u_1v_1)+\eta_0\xi u_1v_1))\\&=&(w_0u_1v_1(2\xi\eta_1-\xi\eta_0-d(\xi)),w_1u_1v_1(2\xi\eta_1-\xi\eta_0-d(\xi))).\hspace{3.9cm}(3.28)\end{eqnarray*}
Therefore, (1.13) holds if and only if
$$2\xi\eta_1=\xi\eta_0+d(\xi).\eqno(3.29)$$

We summarize the above result as the following theorem.
\psp

{\bf Theorem 3.3} {\it Let} $({\cal A},\cdot)$ {\it be a commutative associative algebra and let} $d$ {\it be a derivation of} ${\cal A}$. {\it We set a space} $\td{A}$ {\it as in (3.15) and (3.16)}. {\it For any given three elements} $\xi,\eta_0,\eta_1\in {\cal A}$ {\it satisfying (3.29), we define two algebraic operations} $[\cdot,\cdot]$ {\it and} $\circ$ {\it on} $\td{\cal A}$ {\it as in (3.17) and (3.18). Then the family} $(\td{\cal A},[\cdot,\cdot],\circ)$ {\it forms a Gel'fand-Dorfman super bialgebra}.
\psp

Let $({\cal A},\cdot)$ be a commutative associative algebra and let $d_1,d_2,d_3$ be mutually commutative derivations  of $({\cal A},\cdot)$. Define
$$[u,v]_{i,j}=d_i(u)d_j(v)-d_j(u)d_i(v),\;\;[u,v]_i=ud_i(v)-d_i(u)v\qquad\for\;\;u,v\in{\cal A};\eqno(3.30)$$
$$u\circ_{i,b}v=u(d_i+b)(v)\qquad\for\;\;u,v\in{\cal A},\eqno(3.31)$$
where $b\in\Bbb{F}$ is a constant. It can be verified that $({\cal A},[\cdot,\cdot]_{i,j})$ forms a Lie algebra. Moreover, $({\cal A},\circ_{i,b})$ form Novikov algebras by Corollary 2.6 in [X4]. We want to use these algebras to construct new Gel'fand-Dorfman bialgebras. 

For $u,v,w\in {\cal A}$, we denote 
\begin{eqnarray*}& &[u,v,w]_{i,j,l}=[[u,v]_{i,j},w]_l+[[u,v]_l,w]_{i,j}+[[v,w]_{i,j},u]_l\\& &+[[v,w]_l,u]_{i,j}+[[w,u]_{i,j},v]_l+[[w,u]_l,v]_{i,j}.\hspace{6.7cm}(3.32)\end{eqnarray*}
A tedious calculation shows that 
\begin{eqnarray*}[u,v,w]_{i,j,l}&=&d_i(u)d_l(v)d_j(w)+d_j(u)d_i(v)d_l(w)+d_l(u)d_j(v)d_i(w)\\& &-d_i(u)d_j(v)d_l(w)-d_j(u)d_l(v)d_i(w)-d_l(u)d_i(v)d_j(w)\hspace{2.4cm}(3.33)\end{eqnarray*}
for $u,v,w\in {\cal A}$. By (3.33), we can verify that
$$[\cdot,\cdot,\cdot]_{i,j,j}=0.\eqno(3.34)$$
Thus we have:
\psp

{\bf Proposition 3.4}. {\it The pair} $({\cal A}, [\cdot,\cdot]_{1,2}+\lmd[\cdot,\cdot]_1)$  {\it form a Lie algebra for any} $\lmd\in\Bbb{F}$.
\psp

For $u,v,w\in{\cal A}$, we have:
\begin{eqnarray*}& &[w\circ_{3,b}u,v]_{1,2}-[w\circ_{3,b}v,u]_{1,2}+[w,u]_{1,2}\circ_{3,b}v-[w,v]_{1,2}\circ_{3,b}u-w\circ_{3,b}[u,v]_{1,2}\\&=&d_1(w)d_3(u)d_2(v)+wd_1d_3(u)d_2(v)+bd_1(w)ud_2(v)+bwd_1(u)d_2(v)\\& &-d_2(w)d_3(u)d_1(v)-wd_2d_3(u)d_1(v)-bd_2(w)ud_1(v)-bwd_2(u)d_1(v)\\& &-d_1(w)d_2(u)d_3(v)-wd_2(u)d_1d_3(v)-bd_1(w)d_2(u)v-bwd_2(u)d_1(v)\\&&+d_2(w)d_1(u)d_3(v)+wd_1(u)d_2d_3(v)+bd_2(w)d_1(u)v+bwd_1(u)d_2(v)\\& &+d_1(w)d_2(u)d_3(v)+bd_1(w)d_2(u)v-d_2(w)d_1(u)d_3(v)-bd_2(w)d_1(u)v\\& &-d_1(w)d_3(u)d_2(v)-bd_1(w)ud_2(v)+d_2(w)d_3(u)d_1(v)+bd_2(w)ud_1(v)\\& &-wd_1d_3(u)d_2(v)-wd_1(u)d_2d_3(v)-bwd_1(u)d_2(v)\\& &+wd_2(u)d_1d_3(v)+wd_2d_3(u)d_1(v)+bwd_2(u)d_1(v)\\&=&bw[u,v]_{1,2},\hspace{12cm}(3.35)\end{eqnarray*}
\begin{eqnarray*}& &[w\circ_{2,b}u,v]_1-[w\circ_{2,b}v,u]_1+[w,u]_1\circ_{2,b}v-[w,v]_1\circ_{2,b}u-w\circ_{2,b}[u,v]_1\\&=&wd_2(u)d_1(v)+bwud_1(v)-d_1(w)d_2(u)v-wd_1d_2(u)v-bd_1(w)uv-bwd_1(u)v\\&&-wd_1(u)d_2(v)-bwd_1(u)v+d_1(w)ud_2(v)+wud_1d_2(v)+bd_1(w)uv+bwud_1(v)\\&&+wd_1(u)d_2(v)+bwd_1(u)v-d_1(w)ud_2(v)-bd_1(w)uv-wd_2(u)d_1(v)\hspace{4cm}\end{eqnarray*}
\begin{eqnarray*}& &-bwud_1(v)+d_1(w)d_2(u)v+bd_1(w)uv-bwud_1(v)-wd_2(u)d_1(v)\\& &-wud_1d_2(v)+bwd_1(u)v+wd_1d_2(u)v+wd_1(u)d_2(v)\\&=& w[u,v]_{1,2}.\hspace{12cm}(3.36)\end{eqnarray*}
Proposition 3.4 and the above two expressions imply the following theorem. 
\psp

{\bf Theorem 3.5}. {\it The families} $({\cal A},[\cdot,\cdot]_{1,2}+[\cdot,\cdot]_2,\circ_{2,0})$ {\it and} $({\cal A},[\cdot,\cdot]_{2,1}+b[\cdot,\cdot]_1,\circ_{2,b})$ {\it form  Gel'fand-Dorfman bialgebras}.
\psp

{\bf Remark 3.6}. (a) Theorems 3.2 and 3.5 will be used in next sections of classifications.

(b) For a Lie-Poisson superalgebra $({\cal A},\cdot,[\cdot,\cdot])$, the family $({\cal A},[\cdot,\cdot],\cdot)$ in general does not satisfy (1.13). So it in general does not form a Gel'fand-Dorfman super bialgebra.

\section{Classification I}

In this section, we shall classify the Lie algebras over the simple Novikov algebra $({\cal A}_{\Delta,\{0\}},\circ)$ defined in (1.15).

 For convenience, we redenote
$$x_{\al}=x_{\al,0}\qquad\for\;\;\al\in \Delta;\qquad\; {\cal A}={\cal A}_{\Delta,\{0\}}.\eqno(4.1)$$
Now (1.15) becomes:
$$x_{\al}\circ x_{\be}=(\be+b)x_{\al+\be}\qquad\;\for\;\;\al,\be\in\Delta.\eqno(4.2)$$
Assume that $({\cal A},[\cdot,\cdot])$ is a Lie algebra over $({\cal A},\circ)$. Set
$$[x_{\al},x_{\be}]=\sum_{\sgm\in \Delta}a_{\al,\be}^{\sgm}x_{\sgm}\qquad\for\;\;\al,\be\in\Delta,\eqno(4.3)$$
where $a_{\al,\be}^{\sgm}\in\Bbb{C}$ are the structure constants. The skew-symmetry of Lie algebra implies:
$$a_{\al,\be}^{\sgm}=-a_{\be,\al}^{\sgm}\qquad\;\for\;\;\al,\be,\sgm\in \Delta.\eqno(4.4)$$
For any $\al,\be,\gm\in \Delta$, by (1.11) and (4.4),
\begin{eqnarray*}& &[x_{\gm}\circ x_{\al},x_{\be}]-[x_{\gm}\circ x_{\be},x_{\al}]+[x_{\gm},x_{\al}]\circ x_{\be}-[x_{\gm},x_{\be}]\circ x_{\al}-x_{\gm}\circ [x_{\al},x_{\be}]\hspace{2cm}\end{eqnarray*}
\begin{eqnarray*}&=&\sum_{\sgm\in\Delta}\{[(\al+b)a_{\al+\gm,\be}^{\sgm}-(\be+b)a_{\be+\gm,\al}^{\sgm}]x_{\sgm}+(\be+b)a_{\gm,\al}^{\sgm}x_{\sgm+\be}-(\al+b)a_{\gm,\be}^{\sgm}x_{\sgm+\al}\\& &-(\sgm+b)a_{\al,\be}^{\sgm}x_{\sgm+\gm}\}\\&=&
\sum_{\sgm\in\Delta}[(\al+b)(a_{\al+\gm,\be}^{\sgm}-a_{\gm,\be}^{\sgm-\al})+(\be+b)(a_{\al,\be+\gm}^{\sgm}-a_{\al,\gm}^{\sgm-\be})-(\sgm+b-\gm)a_{\al,\be}^{\sgm-\gm}]x_{\sgm}\\&=&0,\hspace{13.7cm}(4.5)\end{eqnarray*}
which is equivalent to:
$$(\al+b)(a_{\al+\gm,\be}^{\sgm}-a_{\gm,\be}^{\sgm-\al})+(\be+b)(a_{\al,\be+\gm}^{\sgm}-a_{\al,\gm}^{\sgm-\be})-(\sgm+b-\gm)a_{\al,\be}^{\sgm-\gm}=0\eqno(4.6)$$
for $\al,\be,\gm,\sgm\in\Delta$.

Letting $\al=\gm=0$ in (4.6), we have:
$$(\be+b)a_{0,\be}^{\sgm}-(\sgm+b)a_{0,\be}^{\sgm}=0\qquad\for\;\;\be,\sgm\in\Delta,\eqno(4.7)$$
which is equivalent to:
$$(\be-\sgm)a_{0,\be}^{\sgm}=0\qquad\for\;\;\be,\sgm\in\Delta.\eqno(4.8)$$
Thus
$$a_{0,\be}^{\sgm}=0\qquad\for\;\;\be,\sgm\in \Delta;\;\be\neq\sgm.\eqno(4.9)$$
We denote 
$$\vf(\be)=a_{0,\be}^{\be}\qquad\for\;\;\be\in\Delta. \eqno(4.10)$$
Obviously $\vf(0)=0$ by (4.4). By (4.9), 
$$[x_0,x_{\be}]=\vf(\be)x_{\be}\qquad\for\;\;\be\in \Delta.\eqno(4.11)$$

Letting $\al=0,\;\sgm=\be+\gm$ in (4.6), we obtain:
$$(\be+b)(a_{0,\be+\gm}^{\be+\gm}-a_{0,\gm}^{\gm}-a_{0,\be}^{\be})=0,\eqno(4.12)$$
which implies:
$$\vf(\be+\gm)=\vf(\be)+\vf(\gm)\qquad \for\;\;\be,\gm\in \Delta;\;\be\neq -b\;\mbox{or}\;\gm\neq -b.\eqno(4.13)$$
If $b\in\Delta$, we have
$$0=\vf(b+(-b))=\vf(b)+\vf(-b).\eqno(4.14)$$
Hence by (4.13) and (4.14),
$$\vf(-b)+\vf(-b)=\vf(-b)+(\vf(b)+\vf(-2b))=\vf(-2b).\eqno(4.15)$$
By (4.13) and (4.15), $\vf:\Delta\rightarrow\Bbb{C}^+$ is a group homomorphism. 

We define an operator ${\cal D}$ on ${\cal A}$ by:
$${\cal D}(u)= x_0\circ u\qquad\for\;\;u\in{\cal A}.\eqno(4.16)$$
Then
$$\Bbb{C}x_{\al}=\{u\in {\cal A}\mid {\cal D}(u)=(\al+b)u\}.\eqno(4.17)$$
Letting $w=x_0,\;u=x_{\al}$ and $v=x_{\be}$ in (1.11) for $\al,\be\in \Delta$, we get:
$$(\al+b)[x_{\al},x_{\be}]-(\be+b)[x_{\be},x_{\al}]+\vf(\al)(\be+b)x_{\al+\be}-\vf(\be)(\al+b)x_{\be+\al}-{\cal D}[x_{\al},x_{\be}]=0,\eqno(4.18)$$
which is equivalent to:
$$({\cal D}-\al-\be-2b)[x_{\al},x_{\be}]=(\vf(\al)(\be+b)-\vf(\be)(\al+b))x_{\al+\be}.\eqno(4.19)$$

If $b=0$, then (4.17) and (4.19) imply:
$$[x_{\al},x_{\be}]=a_{\al,\be}^{\al+\be}x_{\al+\be},\eqno(4.20)$$
and 
$$\vf(\al)\be=\vf(\be)\al\eqno(4.21)$$
for $\al,\be\in \Delta$. When $\Delta=\{0\}$, the classification is trivial. So we assume $\Delta\neq\{0\}$. Let $\al_0\in\Delta$ be any fixed nonzero element and set
$$a={\vf(\al_0)\over a_0}.\eqno(4.22)$$
By (4.21), we have:
$$\vf(\be)=a\be\qquad\for\;\;\be\in\Delta.\eqno(4.23)$$
Set
$$\phi(\al,\be)=a_{\al,\be}^{\al+\be}+a(\al-\be)\qquad\for\;\;\al,\be\in\Delta.\eqno(4.24)$$
By (4.10) and (4.23),
$$\phi(0,\be)=0\qquad\for\;\;\be\in\Delta.\eqno(4.25)$$
Moreover, (4.4) shows that $\phi$ is skew-symmetric. 
Letting $\sgm=\al+\be+\gm$ in (4.6) for $\al,\be,\gm\in\Delta$, we have:
\begin{eqnarray*}&&\al(\phi(\al+\gm,\be)+a(\be-\al-\gm)-\phi(\gm,\be)-a(\be-\gm))\\& &+\be(\phi(\al,\be+\gm)+a(\gm+\be-\al)-\phi(\al,\gm)-a(\gm-\al))\\& &-(\al+\be)(\phi(\al,\be)+a(\be-\al))=0,\hspace{8cm}(4.26)\end{eqnarray*}
which is equivalent to:
$$\al(\phi(\al+\gm,\be)-\phi(\gm,\be)-\phi(\al,\be))=\be(\phi(\be+\gm,\al)-\phi(\gm,\al)-\phi(\be,\al)).\eqno(4.27)$$
For $\al,\be,\gm\in \Delta$, we set
$$S_0(\al,\be,\gm)=\left\{\begin{array}{ll}0&\mbox{if}\;\;\al=0,\\\al^{-1}(\phi(\be+\gm,\al)-\phi(\gm,\al)-\phi(\be,\al))&\mbox{if}\;\;\al\neq 0.\end{array}\right.\eqno(4.28)$$
Note (4.25) and (4.27) imply that 
$$S_0(\cdot,\cdot,\cdot):\;\Delta\times\Delta\times\Delta \rightarrow \Bbb{C}\;\;\mbox{is s symmetric map}.\eqno(4.29)$$
Furthermore,
$$\phi(\be+\gm,\al)=\phi(\gm,\al)+\phi(\be,\al)+\al S_0(\al,\be,\gm)\qquad\for\;\;\al,\be,\gm\in\Delta.\eqno(4.30)$$
So the map $S_0$ measures the nonlinearity of $\phi$.

By (4.20) and (4.24), we have:
$$[x_{\al},x_{\be}]=(\phi(\al,\be)+a(\be-\al))x_{\al+\be}\qquad\for\;\;\al,\be\in \Delta.\eqno(4.31)$$
Thus for $\al,\be,\gm$, the Jacobi identity of Lie algebra imply:
\begin{eqnarray*}& &[[x_{\al},x_{\be}],x_{\gm}]+[[x_{\be},x_{\gm}],x_{\al}]+[[x_{\gm},x_{\al}],x_{\be}]\\&=&(\phi(\al,\be)+a(\be-\al))(\phi(\al+\be,\gm)+a(\gm-\al-\be))x_{\al+\be+\gm}\\& &+(\phi(\be,\gm)+a(\gm-\be))(\phi(\be+\gm,\al)+a(\al-\be-\gm))x_{\al+\be+\gm}\\&&+(\phi(\gm,\al)+a(\al-\gm))(\phi(\gm+\al,\be)+a(\be-\gm-\al))x_{\al+\be+\gm}\\&=&(\phi(\al,\be)+a(\be-\al))(\phi(\al,\gm)+\phi(\be,\gm)+\gm S_0(\al,\be,\gm)+a(\gm-\al-\be))x_{\al+\be+\gm}\\&&+(\phi(\be,\gm)+a(\gm-\be))(\phi(\be,\al)+ \phi(\gm,\al)+\al S_0(\al,\be,\gm)+a(\al-\be-\gm))x_{\al+\be+\gm}\\&&+(\phi(\gm,\al)+a(\al-\gm))(\phi(\gm,\be)+\phi(\al,\be)+\be S_0(\al,\be,\gm)+a(\be-\gm-\al))x_{\al+\be+\gm}\\&=&\{\phi(\al,\be)\phi(\al,\gm)+\phi(\al,\be)\phi(\be,\gm)+\phi(\be,\gm)\phi(\be,\al)+\phi(\be,\gm)\phi(\gm,\al)\\& &+\phi(\gm,\al)\phi(\gm,\be)+\phi(\gm,\al)\phi(\al,\be)-a(\phi(\al,\be)\gm+\phi(\be,\gm)\al+\phi(\gm,\al)\be)
+[\gm\phi(\al,\be)\\& &+a\gm(\be-\al)+\al\phi(\be,\gm)+a\al(\gm-\be)+\be\phi(\gm,\al)+a\be(\al-\gm)]S_0(\al,\be,\gm)\}x_{\al+\be+\gm}\\&=&(\gm\phi(\al,\be)+\al\phi(\be,\gm)+\be\phi(\gm,\al))(S_0(\al,\be,\gm)-a)x_{\al+\be+\gm}=0,\hspace{2.9cm}(4.32)\end{eqnarray*}
which is equivalent to
$$(\gm\phi(\al,\be)+\al\phi(\be,\gm)+\be\phi(\gm,\al))(S_0(\al,\be,\gm)-a)=0.\eqno(4.33)$$

Note the above arguments are reversible. We summarize the above result as the following theorem.
\psp

{\bf Theorem 4.1}. {\it Any Lie algebra over the simple Novikov algebra} $({\cal A}_{\Delta,\{0\}},\circ)$ {\it with} $b=0$ {\it has its Lie bracket in the form (4.31), where} $a$ {\it is a constant and} $\phi(\cdot,\cdot):\Delta\times\Delta\rightarrow \Bbb{C}$ {\it is a skew-symmetric map such that there exists a symmetric map} $S_0(\cdot,\cdot,\cdot):\Delta\times\Delta\times\Delta\rightarrow \Bbb{C}$ {\it satisfying (4.30), (4.33). Conversely, for any constant} $a\in\Bbb{C}$ {\it and a skew-symmetric map} $\phi(\cdot,\cdot):\Delta\times\Delta\rightarrow \Bbb{C}$ {\it  such that there exists a symmetric map} $S_0(\cdot,\cdot,\cdot):\Delta\times\Delta\times\Delta\rightarrow \Bbb{C}$ {\it satisfying (4.30) and (4.33), the bracket in (4.31) define a Lie algebra over the simple Novikov algebra} $({\cal A}_{\Delta,\{0\}},\circ)$ {\it with} $b=0$. {\it In particular, the following} $\phi$ {\it satisfies our condition: when} $a=0$, $\phi$ {\it is any skew-symmetric} $\Bbb{Z}$-{\it bilinear form; when} $a\neq 0$,
$$\phi(\al,\be)=\al\vf_0(\be)-\be\vf_0(\al)\qquad\mbox{\it for}\;\;\al,\be\in\Delta,\eqno(4.34)$$
{\it where} $\vf_0:\Delta\rightarrow \Bbb{C}^+$ {\it is a group homomorphism}. 
\psp

Next we consider the case $b\not\in \Delta$. Now by (4.17) and (4.19), we have:
$$-b[x_{\al},x_{\be}]=(\vf(\al)(\be+b)-\vf(\be)(\al+b))x_{\al+\be}\qquad\for\;\;\al,\be\in\Delta.\eqno(4.35)$$
Thus
$$[x_{\al},x_{\be}]={1\over b}(\vf(\be)\al-\vf(\al)\be+b(\vf(\be)-\vf(\al)))x_{\al+\be}\qquad\for\;\;\al,\be\in\Delta.\eqno(4.36)$$
\psp

{\bf Theorem 4.2}. {\it A Lie algebra is a Lie algebra over the simple Novikov algebra} $({\cal A}_{\Delta,\{0\}},\circ)$ {\it with} $ b\not \in \Delta$ {\it if and only if  its Lie bracket has the form (4.36), where} $\vf:\Delta\rightarrow \Bbb{C}^+$ {\it is a group homomorphism}.
\psp

{\it Proof}. We only need to prove sufficient part. Let $\vf:\Delta\rightarrow \Bbb{C}^+$ be a group homomorphism and define the operation $[\cdot,\cdot]$ on ${\cal A}={\cal A}_{\Delta,\{0\}}$ by (4.36). Moreover, we define two operations $d_1$ and $d_2$ on ${\cal A}$ by:
$$d_1(x_{\al})={1\over b}\vf(\al)x_{\al},\;\;d_2(x_{\al})=\al x_{\al}\qquad\for\;\;\al\in\Delta.\eqno(4.37)$$
Then $d_1$ and $d_2$ are mutually commutative derivations of ${\cal A}$. By Theorem 3.6, $({\cal A},[\cdot,\cdot],\circ)=({\cal A},[\cdot,\cdot]_{2,1}+b[\cdot,\cdot]_1,\circ_{2,b})$ forms a Gel'fand-Dorfman bialgebra. $\qquad\Box$
\psp 

Finally we consider the case $0\neq b\in \Delta$. Again by (4.17) and (4.19), we have:
$$[x_{\al},x_{\be}]=a_{\al,\be}^{\al+\be+b}x_{\al+\be+b}+a_{\al,\be}^{\al+\be}x_{\al+\be}\eqno(4.38)$$
with
$$a^{\al+\be}_{\al,\be}={1\over b}(\vf(\be)\al-\vf(\al)\be+b(\vf(\be)-\vf(\al))).\eqno(4.39)$$
We set
$$\theta(\al,\be)=a_{\al,\be}^{\al+\be+b}\qquad\for\;\;\al,\be\in\Delta.\eqno(4.40)$$
By (4.4) and (4.9), $\theta(\cdot,\cdot):\Delta\rightarrow \Bbb{C}$ is a skew-symmetric map and
$$\theta(0,\be)=0\qquad\for\;\;\be\in\Delta.\eqno(4.41)$$
Letting $\sgm=\al+\be+\gm+b$ in (4.6) for $\al,\be,\gm\in\Delta$, we obtain:
$$(\al+b)(\theta(\al+\gm,\be)-\theta(\gm,\be)-\theta(\al,\be))=(\be+b)(\theta(\be+\gm,\al)-\theta(\gm,\al)-\theta(\be,\al)).\eqno(4.42)$$
In particular, for $\be=-b$ and $\al\neq -b$ or $\gm\neq -b$, we have:
$$\theta(\al+\gm,-b)=\theta(\gm,-b)+\theta(\al,-b).\eqno(4.43)$$
As (4.13)-(4.15), we can prove that (4.43) holds for any $\al,\gm\in \Delta$. 
We set
$$S_b(\al,\be,\gm)=\left\{\begin{array}{ll}0&\mbox{if}\;\;\al=-b,\\(\al+b)^{-1}(\theta(\be+\gm,\al)-\theta(\gm,\al)-\theta(\be,\al))&\mbox{if}\;\;\al\neq -b\end{array}\right.\eqno(4.44)$$
for $\al,\be,\gm\in \Delta$. Note (4.42) and (4.43) imply that 
$$S_b(\al_1,\al_2,\al_3)\;\;\mbox{is  symmetric with respect to}\;\;\al_i\;\;\mbox{and}\;\;\al_j\eqno(4.45)$$
whenever $(\al_i+b)(\al_j+b)\neq 0$. Moreover,
$$\theta(\be+\gm,\al)=\theta(\gm,\al)+\theta(\be,\al)+(\al+b)S_b(\al,\be,\gm)\qquad\for\;\;\al,\be,\gm\in\Delta.\eqno(4.46)$$

For $\al,\be,\gm$, by (4.45) and (4.46), the Jacobi identity of Lie algebra imply:
\begin{eqnarray*}& &[[x_{\al},x_{\be}],x_{\gm}]+[[x_{\be},x_{\gm}],x_{\al}]+[[x_{\gm},x_{\al}],x_{\be}]\\&=&{1\over b}\{[b\theta(\al,\be)x_{\al+\be+b}+((\al+b)\vf(\be)-(\be+b)\vf(\al))x_{\al+\be},x_{\gm}]\\&&+[b\theta(\be,\gm)x_{\be+\gm+b}+((\be+b)\vf(\gm)-(\gm+b)\vf(\be))x_{\be+\gm},x_{\al}]\\& &+[b\theta(\gm,\al)x_{\gm+\al+b}+((\gm+b)\vf(\al)-(\al+b)\vf(\gm))x_{\al+\gm},x_{\be}]\}\\&=&[\theta(\al,\be)\theta(\al+\be+b,\gm)+\theta(\be,\gm)\theta(\be+\gm+b,\al)+\theta(\gm,\al)\theta(\gm+\al+b,\be)]x_{\al+\be+\gm+2b}\\& &+{1\over b}[\theta(\al,\be)((\al+\be+2b)\vf(\gm)-(\gm+b)\vf(\al+\be+b))+\theta(\al+\be,\gm)((\al+b)\vf(\be)\\& &-(\be+b)\vf(\al))+\theta(\be,\gm)((\be+\gm+2b)\vf(\al)-(\al+b)\vf(\be+\gm+b))\\& &+\theta(\be+\gm,\al)((\be+b)\vf(\gm)-(\gm+b)\vf(\be))+\theta(\gm,\al)((\gm+\al+2b)\vf(\be)\\& &-(\be+b)\vf(\gm+\al+b))+\theta(\gm+\al,\be)((\gm+b)\vf(\al)-(\al+b)\vf(\gm))]x_{\al+\be+\gm+b}\\&=&[\theta(\al,\be)\theta(\al+\be+b,\gm)+\theta(\be,\gm)\theta(\be+\gm+b,\al)+\theta(\gm,\al)\theta(\gm+\al+b,\be)]x_{\al+\be+\gm+2b}\\&&-{1\over b}[(\gm+b)\theta(\al,\be)+(\al+b)\theta(\be,\gm)+(\be+b)\theta(\gm,\al)]\vf(b))x_{\al+\be+\gm+b}\\&=&0,\hspace{13.5cm}(4.47)\end{eqnarray*}
which is equivalent to:
$$\theta(\al,\be)\theta(\al+\be+b,\gm)+\theta(\be,\gm)\theta(\be+\gm+b,\al)+\theta(\gm,\al)\theta(\gm+\al+b,\be)=0\eqno(4.48)$$
and 
$$[(\gm+b)\theta(\al,\be)+(\al+b)\theta(\be,\gm)+(\be+b)\theta(\gm,\al)]\vf(b)=0.\eqno(4.49)$$
When $\gm=0$ in (4.49), we get:
$$\theta(\al,\be)\vf(b)=0\qquad\for\;\;\al,\be\in\Delta.\eqno(4.50)$$
Thus  (4.49) is equivalent to (4.50) and $\theta\equiv 0$ if $\vf(b)\neq 0$. 

If $\vf(b)=0$ and $\theta$ is $\Bbb{Z}$-bilinear, then (4.42) naturally holds and (4.48) is equivalent to
$$\theta(\al,\be)\theta(b,\gm)+\theta(\be,\gm)\theta(b,\al)+\theta(\gm,\al)\theta(b,\be)=0\qquad\for\;\;\al,\be,\gm\in\Delta\eqno(4.51)$$
by the skew-symmetry of $\theta$. If there exists $\gm_0\in\Delta$ such that $\theta(b,\gm_0)\neq 0$, then
$$\theta(\al,\be)=\theta(b,\al){\theta(\gm_0,\be)\over \theta(b,\gm_0)}-\theta(b,\be){\theta(\gm_0,\al)\over \theta(b,\gm_0)}\qquad\for\;\;\for\;\;\al,\be\in \Delta.\eqno(4.52)$$
Thus the solution of (4.51) is:
$$\theta\;\;\mbox{is any skew-symmetric}\;\Bbb{Z}\mbox{-blinear form such that}\;\;b\in\mbox{Rad}_{\theta}\eqno(4.53)$$
or
$$\theta(\al,\be)=\vf_1(\al)\vf_2(\be)-\vf_1(\al)\vf_2(\al)\qquad\for\;\;\al,\be\in\Delta,\eqno(4.54)$$
where $\vf_1,\vf_2:\Delta\rightarrow \Bbb{C}^+$ are group homorphisms such that
$$\vf_1(b)=0,\qquad\vf_2(b)=-1.\eqno(4.55)$$

Since the above arguments are reversible, we have the following theorem:

\psp

{\bf Theorem 4.3}. {\it Any Lie algebra over the simple Novikov algebra} $({\cal A}_{\Delta,\{0\}},\circ)$ {\it with} $0\neq b\in \Delta$ {\it has its Lie bracket as follows:}
$$[x_{\al},x_{\be}]=\theta(\al,\be)x_{\al+\be+b}+{1\over b}((\al+b)\vf(\be)-(\be+b)\vf(\al))x_{\al+\be}\qquad\for\;\;\al,\be\in\Delta.\eqno(4.56)$$
{\it where} $\vf:\Delta\rightarrow \Bbb{C}^+$ {\it is a group homomorphism,} $\theta\equiv 0$ {\it if} $\vf(b)\neq 0$ {\it and} $\theta:\Delta\times\Delta\rightarrow\Bbb{C}$ {\it is a skew-symmetric map satisfying (4.42) and (4.48) if} $\vf(b)=0$. {\it Conversely, for any given group homomorphism} $\vf:\Delta\rightarrow \Bbb{C}^+$ {\it and a skew-symmetric map} $\theta:\Delta\times\Delta\rightarrow\Bbb{C}$ {\it such that} $\theta\equiv 0$ {\it if} $\vf(b)\neq 0$ {\it and (4.42), (4.48) hold if} $\vf(b)=0$, {\it (4.56) defines a Lie algebra over the simple Novikov algebra} $({\cal A}_{\Delta,\{0\}},\circ)$ {\it with} $0\neq b\in \Delta$. {\it In particualr, for any  group homomorphism} $\vf:\Delta\rightarrow \Bbb{C}^+$ {\it such that} $\vf(b)=0$ {\it and a skew-symmetric} $\Bbb{Z}${\it -bilinear map} $\theta:\Delta\times\Delta\rightarrow\Bbb{C}$ {\it in (4.53) or (4.54)}, {\it (4.56) defines a Lie algebra over the simple Novikov algebra} $({\cal A}_{\Delta,\{0\}},\circ)$ {\it with} $0\neq b\in \Delta$.
\psp

{\bf Remark 4.4}. The following Lie algebraic structure on ${\cal A}={\cal A}_{\Delta,\{0\}}$ seems interesting itself, although it not directly related to Gel'fand-Dorfman bialgebras. In fact, it can also be viewed as generalizations of Block algebras (cf. [B]). Assume that $0\neq b\in\Delta$. Let $\phi(\cdot,\cdot):\Delta\times\Delta \rightarrow\Bbb{C}$ be a skew-symmetric $\Bbb{Z}$-bilinear form such that $b\not\in\mbox{Rad}_{\phi}$ and let $\vf:\Delta\rightarrow \Bbb{C}^+$ be a group homomorphism. We have the following Lie bracket on ${\cal A}$:
$$[x_{\al},x_{\be}]=(\vf(\al)\phi(b,\be)-\vf(\be)\phi(b,\al)x_{\al+\be+b}+\phi(\al,\be)x_{\al+\be}\qquad\for\;\;\al,\be\in\Delta.\eqno(4.57)$$

\section{Classification II}

In this section, we shall  classify the Lie algebras over the simple Novikov algebra $({\cal A}_{\Delta,\Bbb{N}},\circ)$ defined in (1.15) with $b\not\in\Delta$. Several families of Lie algebras over the simple Novikov algebra $({\cal A}_{\Delta,\Bbb{N}},\circ)$  with $ b\in\Delta$ will be constructed.

We start with the general case of arbitrary $b$. Let $({\cal A}_{\Delta,\Bbb{N}},[\cdot,\cdot])$ be a Lie algebra over the Novikov algebra $({\cal A}_{\Delta,\Bbb{N}},\circ)$. Set
$$[x_{\al,i},x_{\be,j}]=\sum_{\sgm\in\Delta,\;k\in\Bbb{N}}a_{\al,i;\be,j}^{\sgm,k}x_{\sgm,k}\qquad\for\;\;\al,\be\in\Delta,\;i,j\in\Bbb{N}.\eqno(5.1)$$
Then the skew-symmetry of Lie algebra shows
$$a_{\al,i;\be,j}^{\sgm,k}=-a_{\be,j;\al,i}^{\sgm,k}\qquad\for\;\;\al,\be,\sgm\in\Delta,\;i,j,k\in\Bbb{N}.\eqno(5.2)$$
For $\al,\be,\gm\in\Delta$ and $i,j,l\in\Bbb{N}$, by (1.11) and (1.15),
\begin{eqnarray*}& &[x_{\gm,l}\circ x_{\al,i},x_{\be,j}]-[x_{\gm,l}\circ x_{\be,j}, x_{\al,i}]+[x_{\gm,l},x_{\al,i}]\circ x_{\be,j}\\&&-[x_{\gm,l},x_{\be,j}]\circ x_{\al,i}- x_{\gm,l}\circ [x_{\al,i},x_{\be,j}]\\&=&\sum_{\sgm\in\Delta,\;k\in\Bbb{N}}\{(\al+b)(a_{\al+\gm,i+l;\be,j}^{\sgm,k}x_{\sgm,k}-a_{\gm,l;\be,j}^{\sgm,k}x_{\sgm+\al,k+i})+(\be+b)(a_{\al,i;\be+\gm,j+l}^{\sgm,k}x_{\sgm,k}\\& &+a_{\gm,l;\al,i}^{\sgm,k}x_{\sgm+\be,k+j})+i(a_{\al+\gm,i+l-1;\be,j}^{\sgm,k}x_{\sgm,k}-a_{\gm,l;\be,j}^{\sgm,k}x_{\sgm+\al,k+i-1})+j(a_{\al,i;\be+\gm,j+l-1}^{\sgm,k}x_{\sgm,k}\\& &+a_{\gm,l;\al,i}^{\sgm,k}x_{\sgm+\be,k+j-1})-(\sgm+b)a_{\al,i;\be,j}^{\sgm,k}x_{\sgm+\gm,k+l}-ka_{\al,i;\be,j}^{\sgm,k}x_{\sgm+\gm,k+l-1}\}\hspace{3cm}\end{eqnarray*}
\begin{eqnarray*}&=&\sum_{\sgm\in\Delta,\;k\in\Bbb{N}}\{(\al+b)(a_{\al+\gm,i+l;\be,j}^{\sgm,k}-a_{\gm,l;\be,j}^{\sgm-\al,k-i})+(\be+b)(a_{\al,i;\be+\gm,j+l}^{\sgm,k}+a_{\gm,l;\al,i}^{\sgm-\be,k-j})\\& &+i(a_{\al+\gm,i+l-1;\be,j}^{\sgm,k}-a_{\gm,l;\be,j}^{\sgm-\al,k+1-i})+j(a_{\al,i;\be+\gm,j+l-1}^{\sgm,k}+a_{\gm,l;\al,i}^{\sgm-\be,k+1-j})\\& &-(\sgm-\gm+b)a_{\al,i;\be,j}^{\sgm-\gm,k-l}-(k+1-l)a_{\al,i;\be,j}^{\sgm-\gm,k+1-l}\}x_{\sgm,k}= 0,\hspace{4.3cm}(5.3)\end{eqnarray*}
which is equivalent to:
\begin{eqnarray*}& &(\al+b)(a_{\al+\gm,i+l;\be,j}^{\sgm,k}-a_{\gm,l;\be,j}^{\sgm-\al,k-i})+(\be+b)(a_{\al,i;\be+\gm,j+l}^{\sgm,k}+a_{\gm,l;\al,i}^{\sgm-\be,k-j})\\&&+i(a_{\al+\gm,i+l-1;\be,j}^{\sgm,k}-a_{\gm,l;\be,j}^{\sgm-\al,k+1-i})+j(a_{\al,i;\be+\gm,j+l-1}^{\sgm,k}+a_{\gm,l;\al,i}^{\sgm-\be,k+1-j})\\& &-(\sgm-\gm+b)a_{\al,i;\be,j}^{\sgm-\gm,k-l}-(k+1-l)a_{\al,i;\be,j}^{\sgm-\gm,k+1-l}=0\hspace{5.7cm}(5.4)\end{eqnarray*}
for $\al,\be,\gm,\sgm\in\Delta$ and $i,j,l,k\in\Bbb{N}$. 

Letting $\al=\gm=0$ and $i=l=0$ in (5.4), we have:
$$(\be+b)a_{0,0;\be,j}^{\sgm,k}+ja_{0,0;\be,j-1}^{\sgm,k}-(\sgm+b)a_{0,0;\be,j}^{\sgm,k}-(k+1)a_{0,0;\be,j}^{\sgm,k+1}=0,\eqno(5.5)$$
that is,
$$a_{0,0;\be,j}^{\sgm,k+1}={1\over k+1}[(\be-\sgm)a_{0,0;\be,j}^{\sgm,k}+ja_{0,0;\be,j-1}^{\sgm,k}].\eqno(5.6)$$
Thus by mathematical induction, we have:
$$a_{0,0;\be,j}^{\sgm,k}=\sum_{p=0}^k{1\over (k-p)!}(^j_p)(\be-\sgm)^{k-p}a_{0,0;\be,j-p}^{\sgm,0}\qquad\for\;\;\be\in\Delta,\;j,k\in\Bbb{N}.\eqno(5.7)$$
For any given $\be\in\Delta$ and $j\in \Bbb{N}$, since $[x_{0,0},x_{\be,j}]$ is a finite linear combination of $\{x_{\sgm,k}\mid \sgm\in\Delta,\;k\in\Bbb{N}\}$, there exists integer $K>j$ such that
$$a_{0,0;\be,j}^{\sgm,k}=0\qquad\mbox{when}\;\;k\geq K.\eqno(5.8)$$
Hence (5.7) and (5.8) give the following system:
$$\sum_{p=0}^j{1\over (K+q-p)!}(^j_p)(\be-\sgm)^{K+q-p}a_{0,0;\be,j-p}^{\sgm,0}=0\qquad\for\;\;q=0,1,2,...,j.\eqno(5.9)$$
Assume $\be\neq \sgm$. We view $\{a_{0,0;\be,j}^{\sgm,0},a_{0,0;\be,j-1}^{\sgm,0},...,a_{0,0;\be,0}^{\sgm,0}\}$ as unknowns. The coefficient determinant of the system (5.9) is:
\begin{eqnarray*}& &\left|\begin{array}{cccc}{(\be-\sgm)^K\over K!},&{(\be-\sgm)^{K-1}\over (K-1)!}(^j_1)&\cdots,&{(\be-\sgm)^{K-j}\over (K-j)!}(^j_j)\\ \\{(\be-\sgm)^{K+1}\over (K+1)!},&{(\be-\sgm)^K\over K!}(^j_1)&\cdots,&{(\be-\sgm)^{K+1-j}\over (K+1-j)!}(^j_j)\\\vdots&\vdots&\vdots&\vdots\\{(\be-\sgm)^{K+j}\over (K+j)!},&{(\be-\sgm)^{K+j-1}\over (K+j-1)!}(^j_1)&\cdots,&{(\be-\sgm)^K\over K!}(^j_j)\end{array}\right|\\&=&\left(\prod_{p=0}^j(^j_p){(\be-\sgm)^K\over (K+p)!}\right)\left|\begin{array}{cccc}1,&K,&\cdots,& K(K-1)\cdots(K+1-j)\\ 1,&(K+1),&\cdots,& (K+1)K\cdots(K+2-j)\\\vdots&\vdots&\vdots&\vdots\\ 1,&(K+j),&\cdots,& (K+j)(K+j-1)\cdots(K+1)\end{array}\right|\hspace{2cm}\end{eqnarray*}
\begin{eqnarray*}&=&\left[\left(\prod_{p=0}^j(^j_p){(\be-\sgm)^K\over (K+p)!}\left({d\over dz_p}\right)^p\right)\left|\begin{array}{cccc}1,&z_1^K,&\cdots,&z_j^K\\1,&z_1^{K+1},&\cdots,&z_j^{K+1}\\\vdots&\vdots&\vdots&\vdots\\1,&z_1^{K+j},&\cdots,&z_j^{K+j}\end{array}\right|\right]_{z_1=z_2=\cdots z_j=1}\\&=&\left[\left(\prod_{p=0}^j(^j_p){(\be-\sgm)^K\over (K+p)!}\left({d\over dz_p}\right)^p\right)(z_1z_2\cdots z_j)^K\left|\begin{array}{cccc}1,&1,&\cdots,&1\\1,&z_1,&\cdots,&z_j
\\\vdots&\vdots&\vdots&\vdots\\1,&z_1^j,&\cdots,&z_j^j\end{array}\right|\right]_{z_1=z_2=\cdots z_j=1}\\&=&
\left[\left(\prod_{p=0}^j(^j_p){(\be-\sgm)^K\over (K+p)!}\left({d\over dz_p}\right)^p\right)(z_1z_2\cdots z_j)^K(\prod_{i=1}^j(1-z_i))(\prod_{1\leq s<t\leq j}(z_s-z_t))\right]_{z_1=z_2=\cdots z_j=1}\\&=&(-1)^{j(j+1)/2}\prod_{p=0}^j(^j_p){(\be-\sgm)^K\over (K+p)!}\neq 0.\hspace{8.2cm}(5.10)\end{eqnarray*}
Thus we have:
$$a_{0,0;\be,j}^{\sgm,0}=a_{0,0;\be,j-1}^{\sgm,0}=\cdots=a_{0,0;\be,0}^{\sgm,0}=0.\eqno(5.11)$$
By (5.7), we obtain:
$$a_{0,0;\be,j}^{\sgm,k}=0\qquad\for\;\;\be\neq \sgm\in \Delta,\;k\in\Bbb{Z}\eqno(5.12)$$
and
$$a_{0,0;\be,j}^{\be,k}=\left\{\begin{array}{ll}0&\mbox{if}\;\;k>j,\\(^j_k)a_{0,0;\be,j-k}^{\be,0}&\mbox{if}\;\;k\leq j\end{array}\right.\eqno(5.13)$$

Letting $\al=\gm=0,\;\sgm=\be,\;i=0$ and $k=l=1$ in (5.4), we get:
$$(\be+b)(a_{0,0;\be,j+1}^{\be,1}+a_{0,1;0,0}^{0,1-j})+j(a_{0,0;\be,j}^{\be,1}+a_{0,1;0,0}^{0,2-j})-(\be+b)a_{0,0;\be,j}^{\be,0}-a_{0,0;\be,j}^{\be,1}=0.\eqno(5.14)$$
By (5.13) and (5.14),
$$(\be+b)(ja_{0,0;\be,j}^{\be,0}+\delta_{1,j}a_{0,1;0,0}^{0,0})+(j-1)ja_{0,0;\be,j-1}^{\be,0}+j\delta_{2,j}a_{0,1;0,0}^{0,0}=0.\eqno(5.15)$$
When $j=1$ in (5.15), we have:
$$(\be+b)(a_{0,0;\be,1}^{\be,0}+a_{0,1;0,0}^{0,0})=0.\eqno(5.16)$$
Thus
$$a_{0,0;\be,1}^{\be,0}=a_{0,0;0,1}^{0,0}\qquad\for\;\;-b\neq \be\in\Delta\eqno(5.17)$$
by (5.2). When $j=2$ in (5.15),
$$2(\be+b)a_{0,0;\be,2}^{\be,0}+2a_{0,0;\be,1}^{\be,0}+2a_{0,1;0,0}^{0,0}=0,\eqno(5.18)$$
which is equivalent to:
$$a_{0,0;\be,2}^{\be,0}=0\qquad\for\;\;-b\neq \be\in\Delta\eqno(5.20)$$
by (5.2) and (5.17). When $j>2$ in (5.15), we get:
$$j(\be+b)a_{0,0;\be,j}^{\be,0}+(j-1)ja_{0,0;\be,j-1}^{\be,0}=0,\eqno(5.21)$$
which is equivalent to
$$(\be+b)a_{0,0;\be,j}^{\be,0}=(1-j)a_{0,0;\be,j-1}^{\be,0}.\eqno(5.22)$$
Moreover, by (5.20), (5.22) and mathematical induction on $j$, we can prove:
$$a_{0,0;\be,j}^{\be,0}=0\qquad\for\;\;-b\neq \be\in\Delta,\;2\leq j\in\Bbb{N}.\eqno(5.23)$$
When $\be=-b$ in (5.15), we have:
$$(j-1)ja_{0,0;-b,j-1}^{-b,0}+j\delta_{2,j}a_{0,1;0,0}^{0,0}=0,\eqno(5.24)$$
which implies:
$$a_{0,0;-b,1}^{-b,0}=a_{0,0;0,1}^{0,0},\qquad a_{0,0;-b,j}^{-b,0}=0\qquad\for\;\;2\leq j\in\Bbb{N}.\eqno(5.25)$$
Therefore, we have:
$$a_{0,0;\be,1}^{\be,0}=a_{0,0;0,1}^{0,0},\qquad a_{0,0;\be,j}^{\be,0}=0\qquad\for\;\;\be\in\Delta,\;2\leq j\in\Bbb{N}.\eqno(5.26)$$

We denote 
$$a_{0,0;0,1}^{0,0}=\lmd,\;\;a_{0,0;\be,0}^{\be,0}=\vf(\be)\qquad\for\;\;\be\in\Delta.\eqno(5.27)$$
By (5.13) and (5.26),
$$a_{0,0;\be,j}^{\be,j-1}=\lmd j\qquad\for\;\;\be\in\Delta,\;j\in\Bbb{N}.\eqno(5.28)$$
Thus we have:
$$[x_{0,0},x_{\be,j}]=\vf(\be)x_{\be,j}+\lmd jx_{\be,j-1}\qquad\for\;\;\be\in\Delta,\;j\in\Bbb{N}.\eqno(5.29)$$
Note that letting $\al=0,\;\sgm=\be+\gm$ and $i=j=l=k=0$ in (5.4), we get the same equation as (4.12) with $a_{0,\rho}^{\rho}$ replaced by $a_{0,0;\rho,0}^{\rho,0}$ for $\rho\in\Delta$ because $a_{0,0;\be,0}^{\be,1}=0$ by (5.13). So $\vf:\Delta\rightarrow\Bbb{C}^+$ is a group homomorphism by (4.12)-(4.15).

Define an operator ${\cal D}$ on ${\cal A}_{\Delta,\Bbb{N}}$ by
$${\cal D}(u)=x_{0,0}\circ u\qquad\for\;\;u\in{\cal A}_{\Delta,\Bbb{N}}.\eqno(5.30)$$
Then
$$\sum_{p=0}^{j-1}\Bbb{C}x_{\be,p}=\{u\in {\cal A}_{\Delta,\Bbb{N}}\mid ({\cal D}-\be-b)^j(u)=0\}\qquad\for\;\;\be\in\Delta,\;j\in\Bbb{N}.\eqno(5.31)$$
For $\al,\be\in\Delta$ and $i,j\in\Bbb{N}$, letting $w=x_{0,0},\;u=x_{\al,i}$ and $v=x_{\be,j}$ in (1.11), we get:
\begin{eqnarray*}& &[(\al+b)x_{\al,i}+ix_{\al,i-1},x_{\be,j}]-[(\be+b)x_{\be,j}+jx_{\be,j-1},x_{\al,i}]+(\vf(\al)x_{\al,i}+\lmd i x_{\al,i-1})\circ x_{\be,j}\\&&-(\vf(\be)x_{\be,j}+\lmd j x_{\be,j-1})\circ x_{\al,i}-x_{0,0}\circ [x_{\al,i},x_{\be,j}]=0,\hspace{4.9cm}(5.32)\end{eqnarray*}
which is equivalent to:
\begin{eqnarray*} & &({\cal D}-\al-\be-2b)([x_{\al,i},x_{\be,j}])\\&=&i[x_{\al,i-1},x_{\be,j}]+j[x_{\al,i},x_{\be,j-1}]+(\vf(\al)(\be+b)-\vf(\be)(\al+b))x_{\al+\be,i+j}\\& &+[i(\lmd(\be+b)-\vf(\be))+j(\vf(\al)-\lmd(\al+b))]x_{\al+\be,i+j-1}.\hspace{4.2cm}(5.33)\end{eqnarray*}

Assume that $b\neq 0$. By (5.31), (5.33) and mathematical induction on $i+j$, we can prove that
$$[x_{\al,i},x_{\be,j}]=\sum_{k=0}^{i+j}(a_{\al,i;\be,j}^{\al+\be+b,k}x_{\al+\be+b,k}+a_{\al,i;\be,j}^{\al+\be,k}x_{\al+\be,k})\qquad\for\;\;\al,\be\in\Delta,\;i,j\in\Bbb{N},\eqno(5.34)$$
where we treat $x_{\al+\be+b,k}=0$ if $b\not\in\Delta$. 
Note that
$$({\cal D}-\al-\be-2b)(\sum_{k=0}^{i+j}a_{\al,i;\be,j}^{\al+\be,k}x_{\al+\be,k})=\sum_{k=0}^{i+j}a_{\al,i;\be,j}^{\al+\be,k}(-bx_{\al+\be,k}+kx_{\al+\be,k-1}).\eqno(5.35)$$
Hence by (5.33)-(5.35) and mathematical induction on $i+j$, we can prove that
$$a_{\al,i;\be,j}^{\al+\be,i+j}={1\over b}((\al+b)\vf(\be)-(\be+b)\vf(\al))\qquad\for\;\;\al,\be\in\Delta,\;i,j\in\Bbb{N}.\eqno(5.36)$$
Moreover, by (5.34) and (5.35), the coefficients of $x_{\al+\be,i+j-1}$ in (5.33) shows that
\begin{eqnarray*}& &-ba_{\al,i;\be,j}^{\al+\be,i+j-1}+{i+j\over b}((\al+b)\vf(\be)-(\be+b)\vf(\al))\\&=&{i+j\over b}((\al+b)\vf(\be)-(\be+b)\vf(\al))+i(\lmd(\be+b)-\vf(\be))\\& &+j(\vf(\al)-\lmd(\al+b)),\hspace{10cm}(5.37)\end{eqnarray*}
which is equivalent to:
$$a_{\al,i;\be,j}^{\al+\be,i+j-1}={1\over b}[i(\vf(\be)-\lmd(\be+b))+j(\lmd(\al+b)-\vf(\al))]\eqno(5.38)$$
for $\al,\be\in\Delta,\;i,j\in\Bbb{N}$. From the coefficients of $x_{\al+\be,i+j-2}$ in (5.33), we get
\begin{eqnarray*}& &-b a_{\al,i;\be,j}^{\al+\be,i+j-2}+{i+j-1\over b}[i(\vf(\be)-\lmd(\be+b))+j(\lmd(\al+b)-\vf(\al))]\\&=&{i\over b}[(i-1)(\vf(\be)-\lmd(\be+b))+j(\lmd(\al+b)-\vf(\al))]\\& &+{j\over b}[i(\vf(\be)-\lmd(\be+b))+(j-1)(\lmd(\al+b)-\vf(\al))]\\&=&{i+j-1\over b}[i(\vf(\be)-\lmd(\be+b))+j(\lmd(\al+b)-\vf(\al))],\hspace{4.5cm}(5.39)\end{eqnarray*}
which is equivalent to:
$$a_{\al,i;\be,j}^{\al+\be,i+j-2}=0\qquad\for\;\;\al,\be\in\Delta,\;i,j\in\Bbb{N}.\eqno(5.40)$$
Thus by (5.33)-(5.35) and mathematical induction, we can prove that
$$a_{\al,i;\be,j}^{\al+\be,k}=0\qquad\for\;\;\al,\be\in\Delta,\;i,j\in\Bbb{N},\;k\leq i+j-2.\eqno(5.41)$$
\psp

{\bf Theorem 5.1}. {\it A Lie algebra is an algebra over the simple Novikov algebra} $({\cal A}_{\Delta,\Bbb{N}},\circ)$ {\it (cf. (1.15)) with} $b\not\in\Delta$ {\it if and only if its Lie bracket has the following form:}
\begin{eqnarray*}[x_{\al,i},x_{\be,j}]&=&{1\over b}((\al+b)\vf(\be)-(\be+b)\vf(\al))x_{\al+\be,i+j}\\& &+{1\over b}[i(\vf(\be)-\lmd(\be+b))+j(\lmd(\al+b)-\vf(\al))]x_{\al+\be,i+j-1}\hspace{2.3cm}(5.42)\end{eqnarray*}
{\it for} $\al,\be\in\Delta,\;i,j\in\Bbb{N}$, {\it where} $\vf:\Delta\rightarrow \Bbb{C}^+$ {\it is a group homomorphism and} $\lmd\in\Bbb{C}$ {\it is a constant}.

{\it Proof}. The necessity has been proved in the above. 

To prove the sufficiency, we define operators $d_1,d_2$ on ${\cal A}_{\Delta,\Bbb{N}}$ by
$$d_1(u_{\al,i})=\vf(\al)u_{\al,i}+\lmd iu_{\al,i-1},\;\;d_2(u_{\al,i})=\al u_{\al,i}+iu_{\al,i-1}\qquad\for\;\;\al\in\Delta,\;i\in\Bbb{N}.\eqno(5.43)$$
Then $d_1$ and $d_2$ are mutually commutative derivations of $({\cal A}_{\Delta,\Bbb{N}},\cdot)$. Moreover, the algebraic operations defined in (5.42) and (1.15) have the property:
$$b [\cdot,\cdot]=[\cdot,\cdot]_{2,1}+b[\cdot,\cdot]_1,\;\;\circ=\circ_b\eqno(5.44)$$
in terms the notions in (3.30) and (3.31). Thus (1.15) and (5.42) define a Gel'fand-Dorfman bialgebra by Theorem 3.5.$\qquad\Box$
\psp

{\bf Remark 5.2}. (a) Up to this stage, we have not completely classified all the Lie algebras over the simple Novikov algebra $({\cal A}_{\Delta,\Bbb{N}},\circ)$ with $b\in\Delta$. The best informations that we have obtained are as follows. When $0\neq b\in\Delta$,
$$a_{\al,i;\be,j}^{\al+\be+b,k}=\sum_{k=0}^k(^i_p)(^{\;\;\;j}_{k-p})a_{\al,i-p;\be,j+p-k}^{\al+\be+b,0}\qquad\for\;\;\al,\be\in\Delta,\;i,j,k\in\Bbb{N},\eqno(5.45)$$
$$a_{0,1;0,j}^{b,0}={(-1)^j(j+1)!\over 6(2b)^{j-2}}a_{0,1;0,2}^{b,0},\;\;a_{0,i;0,j}^{b,0}={(-1)^{i+j}(j-i)(i+j-1)!\over 6(2b)^{i+j-2}}a_{0,2;0,1}^{b,0}\eqno(5.46)$$
for $1<i,j\in\Bbb{N}$. When $b=0$, for $\al,\be\in\Delta$ and $i,j\in\Bbb{N}$,
$$ a_{\al,i;\be,j}^{\al+\be,k}=0\qquad\for\;\;i+j+1<k\in\Bbb{N},\eqno(5.47)$$
$$a_{\al,i;\be,j}^{\al+\be,i+j+1}=\vf(\al)\be-\vf(\be)\al,\;\;a_{\al,i;\be,j}^{\al+\be,i+j}=a_{\al,0;\be,0}^{\al+\be,0}+i(\lmd \be-\vf(\be))+j(\vf(\al)-\lmd\al),\eqno(5.48)$$
$$a_{\al,i;\be,j}^{\al+\be,k}=\sum_{k=0}^k(^i_p)(^{\;\;\;j}_{k-p})a_{\al,i-p;\be,j+p-k}^{\al+\be,0}\qquad\for\;\;i+j>k\in\Bbb{N},\eqno(5.49)$$
$$[x_{0,i},x_{0,j}]=\lmd(j-i).\eqno(5.50)$$
Equation (5.50) was proved by Osborn and Zelmanov [OZ].

(b) The followings are Lie algebras over $({\cal A}_{\Delta,\Bbb{N}},\circ)$ with $b=0$:

(1) For a skew-symmetric $\Bbb{Z}$-bilinear form $\phi(\cdot,\cdot):\Delta\times\Delta\rightarrow \Bbb{C}$ and a group homomorphism $\vf:\Delta\rightarrow\Bbb{C}^+$, the Lie bracket is defined by:
$$[x_{\al,i},x_{\be,j}]=\phi(\al,\be)x_{\al+\be,i+j}+(i\vf(\be)-j\vf(\al))x_{\al+\be,i+j-1}\eqno(5.51)$$
for $\al,\be\in\Delta,\;i,j\in\Bbb{N}$. This structure is obtained by Theorem 3.8 when $[\cdot,\cdot]=[\cdot,\cdot]_{\phi,2,3}$ with $d_i$ defined in (5.43). 

(2) For a group homomorphism $\vf:\Delta\rightarrow\Bbb{C}^+$ and a nonzero constant $\lmd\in\Bbb{C}$, the Lie bracket is defined as follows:
$$[x_{\al,i},x_{\be,j}]=(\al\vf(\be)-\be\vf(\al))x_{\al+\be,i+j}+[i(\vf(\be)-\lmd\be)+j(\lmd\al-\vf(\al))]x_{\al+\be,i+j-1}\eqno(5.52)$$
for $\al,\be\in\Delta,\;i,j\in\Bbb{N}$. This structure is obtained by Corollary 3.7 with $d_i$ defined in (5.43) and $b=0$. 

(3)  For a group homomorphism $\vf:\Delta\rightarrow\Bbb{C}^+$, the Lie bracket is defined as follows:
$$[x_{\al,i},x_{\be,j}]=(\al\vf(\be)-\be\vf(\al)+\be-\al)x_{\al+\be,i+j}+(i\vf(\be)-j\vf(\al)+j-i)x_{\al+\be,i+j-1}\eqno(5.53)$$
for $\al,\be\in\Delta,\;i,j\in\Bbb{N}$. This structure is obtained by Theorem 3.6. 

(c)  Lie algebras over $({\cal A}_{\Delta,\Bbb{N}},\circ)$ with $0\neq b\in\Delta$:

(1) For a group homomorphism $\vf:\Delta\rightarrow\Bbb{C}^+$ and a constant $\lmd\in\Bbb{C}$, the Lie bracket is defined by (5.42).

(2)  For group homomorphisms $\vf,\vf_1,\vf_2:\Delta\rightarrow\Bbb{C}^+$ such that $\vf(b)=\vf_1(b)=0$
and a constant $\lmd\in\Bbb{C}$, the Lie bracket is defined below:
\begin{eqnarray*}[x_{\al,i},x_{\be,j}]&=&(\vf_1(\al)\vf_2(\be)-\vf_1(\be)\vf_2(\al))x_{\al+\be+b,i+j}\\& &+((\al+b)\vf(\be)-(\be+b)\vf(\al))x_{\al+\be,i+j}\\& &+[i(\vf(\be)-\lmd(\be+b))+j(\lmd(\al+b)-\vf(\al))]x_{\al+\be,i+j-1}\hspace{2.4cm}(5.54)\end{eqnarray*}
for $\al,\be\in\Delta,\;i,j\in\Bbb{N}$. This family of Lie algebras are motivated by Theorem 4.3 and (5.42).

(d) The following Lie algebraic structure on ${\cal A}_{\Delta,\Bbb{N}}$ seems interesting itself, although it is not directly related to Gel'fand-Dorfman bialgebras.
Assume that $0\neq b\in\Delta$. Let $\phi(\cdot,\cdot):\Delta\times\Delta \rightarrow\Bbb{C}$ be a skew-symmetric $\Bbb{Z}$-bilinear form such that $b\not\in\mbox{Rad}_{\phi}$ and let $\vf,\vf_1:\Delta\rightarrow \Bbb{C}^+$ be group homomorphisms such that $\vf(b)=0$. We have the following Lie bracket on ${\cal A}_{\Delta,\Bbb{N}}$:
\begin{eqnarray*}[x_{\al},x_{\be}]&=&(\vf_1(\al)\phi(b,\be)-\vf_1(\be)\phi(b,\be))x_{\al+\be+b}+\phi(\al,\be)x_{\al+\be}\\& &+(i\vf(\be)-j\vf(\al))x_{\al+\be,i+j-1}\hspace{7.6cm}(5.55)\end{eqnarray*}
for $\al,\be\in\Delta,\;i,j\in\Bbb{N}$.

\section{Classification III}

Recall we have a commutative associative algebraic structure on ${\cal A}_{\Delta,\G}$ defined in (1.16).
As indicated in (1.18), the  commutator algebras of all simple Novikov algebras $({\cal A}_{\Delta,\G},\dmd_{\xi})$ for $\xi\in {\cal A}_{\Delta,\G}$ are the same. We ask whether $\dmd_{\xi}$ are the all Novikov algebraic structures on $\xi\in {\cal A}_{\Delta,\G}$ whose commutator algebras are the following Lie algebra:
$$[x_{\al,i},x_{\be,j}]=(\be-\al)x_{\al+\be,i+j}+(j-i)x_{\al+\be,i+j-1}\qquad\for\;\;\al,\be\in\Delta,\;i,j\in\Bbb{N}.\eqno(6.1)$$
We shall give a confirmative answer to the case $\G=\{0\}$.

Let $({\cal A}_{\Delta,\G},\circ)$ be a Novikov algebra whose commutative algebra is given by (6.1). Then $({\cal A}_{\Delta,\G},[\cdot,\cdot],\circ)$ forms Gel'fand-Dorfman bialgebra, that is, (1.11) holds (cf. Theorem 2.3). Set
$$\xi=x_{0,0}\circ x_{0,0}\eqno(6.2)$$
and define another algebraic operation  $\circ_0$ on ${\cal A}_{\Delta,\G}$ by:
$$u\circ_0v=\xi uv\qquad\for\;\;u,v\in {\cal A}_{\Delta,\G}.\eqno(6.3)$$
\psp

{\bf Lemma 6.1}. {\it The family} $({\cal A}_{\Delta,\G},[\cdot,\cdot],\circ_0)$ {\it forms a Gel'fand Dorfman bialgebra}.
\psp

{\it Proof}. It is easily seen that $({\cal A}_{\Delta,\G},\circ_0)$ forms a Novikov algebra. Moreover, (1.11) holds because $\circ_0=\dmd_{\xi}-\dmd_0$ and both $([\cdot,\cdot],\dmd_{\xi})$ and $([\cdot,\cdot],\dmd_0)$ satisfy (1.11) (cf. (1.18) and Theorem 2.3).$\qquad\Box$
\psp

Now we let
$$\star=\circ-\circ_0.\eqno(6.4)$$
Note that the commutator algebra of $({\cal A}_{\Delta,\G},[\cdot,\cdot],\circ_0)$ is a trivial (abelian) Lie algebra. Thus the commutator algebra of the algebra $({\cal A}_{\Delta,\G},\star)$ is also
$({\cal A}_{\Delta,\G},[\cdot,\cdot])$, and $([\cdot,\cdot],\star)$ satisfies (1.11). Specifically, we have:
$$x_{\al,i}\star x_{\be,j}-x_{\be,j}\star x_{\al,i}=[x_{\al,i},x_{\be,j}]=(\be-\al)x_{\al+\be,i+j}+(j-i)x_{\al+\be,i+j-1}\eqno(6.5)$$
for $\al,\be\in\Delta,\;i,j\in\Bbb{N},$
$$[w\star u,v]-[w\star v,u]+[w,u]\star v-[w,v]\star u-w\star [u,v]=0\qquad\for\;\;u,v,w\in{\cal A}_{\Delta,\G},\eqno(6.6)$$
$$x_{0,0}\star x_{0,0}=0.\eqno(6.7)$$
Below we assume that $\G=\{0\}$. We shall determine $({\cal A}_{\Delta,\{0\}},\circ)$ through $({\cal A}_{\Delta,\{0\}},\star)$.

Recall the notion in (4.1). We write
$$x_{\al}\star x_{\be}=\sum_{\sgm\in\Delta} c_{\al,\be}^{\sgm}x_{\sgm}\qquad\for\;\;\al,\be\in\Delta.\eqno(6.8)$$
Then by (6.3) and (6.7), we have:
$$c_{\al,\be}^{\sgm}-c_{\be,\al}^{\sgm}=\dlt_{\sgm,\al+\be}(\be-\al)\qquad\for\;\;\al,\be,\sgm\in\Delta,\eqno(6.9)$$
$$c_{0,0}^{\sgm}=0\qquad\for\;\;\sgm\in\Delta.\eqno(6.10)$$
Moreover, for $\al,\be,\gm\in \Delta$, (6.6) gives
\begin{eqnarray*}& &[x_{\gm}\star x_{\al},x_{\be}]-[x_{\gm}\star x_{\be}, x_{\al}]+[x_{\gm},x_{\al}]\star x_{\be}-[x_{\gm},x_{\be}]\star x_{\al}-x_{\gm}\star [x_{\al},x_{\be}]\\&=&\sum_{\sgm\in\Delta}[(\be-\sgm)c_{\gm,\al}^{\sgm}x_{\sgm+\be}-(\al-\sgm)c_{\gm,\be}^{\sgm}x_{\sgm+\al}+((\al-\gm)c_{\gm+\al,\be}^{\sgm}\\& &-(\be-\gm)c_{\gm+\be,\al}^{\sgm}-(\be-\al)c_{\gm,\al+\be}^{\sgm})x_{\sgm}]\\&=&\sum_{\sgm\in\Delta}[(2\be-\sgm)c_{\gm,\al}^{\sgm-\be}-(2\al-\sgm)c_{\gm,\be}^{\sgm-\al}+(\al-\gm)c_{\gm+\al,\be}^{\sgm}\\& &-(\be-\gm)c_{\gm+\be,\al}^{\sgm}-(\be-\al)c_{\gm,\al+\be}^{\sgm}]x_{\sgm}=0,\hspace{6.7cm}(6.11)\end{eqnarray*}
which is equivalent to:
$$(2\be-\sgm)c_{\gm,\al}^{\sgm-\be}-(2\al-\sgm)c_{\gm,\be}^{\sgm-\al}+(\al-\gm)c_{\gm+\al,\be}^{\sgm}-(\be-\gm)c_{\gm+\be,\al}^{\sgm}-(\be-\al)c_{\gm,\al+\be}^{\sgm}=0\eqno(6.12)$$
for $\al,\be,\gm,\sgm\in\Delta$.

Letting $\al=\gm=0$ and $\be\neq 0$ in (6.12), we get:
$$\sgm c_{0,\be}^{\sgm}-\be c_{\be,0}^{\sgm}-\be c_{0,\be}^{\sgm}=0\eqno(6.13)$$
by (6.10), which implies
$$(\sgm-2\be)c_{0,\be}^{\sgm}=-\be^2\dlt_{\sgm,\be}\eqno(6.14)$$
by (6.9). Hence
$$c_{0,\be}^{\sgm}=0\qquad\for\;\;\be,\sgm\in\Delta,\;\sgm\neq \be,2\be,\eqno(6.15)$$
$$c_{0,\be}^{\be}=\be\qquad\for\;\;\be\in\Delta.\eqno(6.16)$$

Next we let $\gm=0,\;\sgm=2(\al+\be),\;\al\neq \pm\be$ and $\al\be\neq 0$ in (6.12), and obtain:
$$\al c_{\al,\be}^{2(\al+\be)}-\be c_{\be,\al}^{2(\al+\be)}-(\be-\al)c_{0,\al+\be}^{2(\al+\be)}=0,\eqno(6.17)$$
by (6.15). Moreover, by (6.9), (6.17) can be written as:
$$(\al-\be)c_{\al,\be}^{2(\al+\be)}-(\be-\al)c_{0,\al+\be}^{2(\al+\be)}=0,\eqno(6.18)$$
which implies:
$$c_{\al,\be}^{2(\al+\be)}=-c_{0,\al+\be}^{2(\al+\be)}\qquad\for\;\;0\neq \al,\be\in\Delta,\;\;\al\neq \pm \be.\eqno(6.19)$$

When $\al=0,\;\sgm=2(\be+\gm)$ and $\be\gm\neq 0,\;\be\neq \pm \gm$ in (6.12), we have:
$$2(\be+\gm)c_{\gm,\be}^{2(\be+\gm)}-\gm c_{\gm,\be}^{2(\be+\gm)}-(\be-\gm)c_{\gm+\be,0}^{2(\be+\gm)}-\be c_{\gm,\be}^{2(\be+\gm)}=0\eqno(6.20)$$
by (6.15). Moreover, by (6.9) and (6.19), (6.20) implies:
$$-2\be c_{0,\be+\gm}^{2(\be+\gm)}=0.\eqno(6.21)$$
Thus
$$c_{\gm,\be}^{2(\be+\gm)}=c_{0,\be+\gm}^{2(\be+\gm)}=0.\eqno(6.22)$$
For any $0\neq \tau\in\Delta$, we let $\be=2\tau$ and $\gm=-\tau$. Then $\be$ and $\gm$ satisfy our assumption. So
$$c_{0,\tau}^{2\tau}=c_{0,2\tau+(-\tau)}^{2(2\tau+(-\tau))}=0\qquad\for\;\;\tau\in\Delta.\eqno(6.23)$$

Assuming that $\al=0$ and $\sgm\neq \be+\gm$ in (6.12), we get:
$$\sgm c_{\gm,\be}^{\sgm}-\gm c_{\gm,\be}^{\sgm}-\be c_{\gm,\be}^{\sgm}=0\eqno(6.24)$$
by (6.15) and (6.22)-(6.23). Thus
$$c_{\gm,\be}^{\sgm}=0\qquad\for\;\;\be,\gm,\sgm\in\Delta,\;\sgm\neq \be+\gm.\eqno(6.25)$$

Supposing that $\gm=0,\;\al\neq \be$ and $\sgm=\al+\be$ in (6.12), we get:
$$(\be-\al) c_{0,\al}^{\al}+(\be-\al)c_{0,\be}^{\be}+\al c_{\al,\be}^{\al+\be}-\be c_{\be,\al}^{\al+\be}-(\be-\al)c_{0,\al+\be}^{\al+\be}=0,\eqno(6.26)$$
which is equivalent to:
$$(\al-\be)c_{\al,\be}^{\al+\be}=\be(\al-\be).\eqno(6.27)$$
So we have:
$$c_{\al,\be}^{\al+\be}=\be\qquad\for\;\;\al,\be\in\Delta,\;\al\neq\be.\eqno(6.28)$$

Letting $\al=3\be,\;\gm=\be\neq 0$ and $\sgm=5\be$ in (6.12), we get:
$$-3\be c_{\be,3\be}^{4\be}-\be c_{\be,\be}^{2\be}+2\be c_{4\be,\be}^{5\be}+2\be c_{\be,4\be}^{5\be}=0,\eqno(6.29)$$
which is equivalent to:
$$-9\be^2-\be c_{\be,\be}^{2\be}+2\be^2+8\be^2=0\eqno(6.30)$$
by (6.28). Thus we have:
$$c_{\be,\be}^{2\be}=\be\qquad\for\;\;\be\in\Delta.\eqno(6.31)$$
By (6.25), (6.28) and (6.31), we have:
$$x_{\al}\star x_{\be}=\be x_{\al+\be}\qquad\for\;\;\al,\be\in \Delta.\eqno(6.32)$$
Hence by (6.2)-(6.4), Lemma 6.1 and (6.32), we obtain the main result in this section:
\psp

{\bf Theorem 6.2}. {The set} $\{({\cal A}_{\Delta,\{0\}},\dmd_{\xi})\mid\xi\in{\cal A}_{\Delta,\{0\}}\}$ {\it enumerates all the Novikov algebraic structures over} ${\cal A}_{\Delta,\{0\}}$ {\it whose commutator algebras are given by (6.1)}.
\vspace{1cm}

\noindent{\Large \bf References}

\hspace{0.5cm}

\begin{description}

\item[{[BN]}] A. A. Balinskii and S. P. Novikov, Poisson brackets of hydrodynamic type, Frobenius algebras and Lie algebras, {\it Soviet Math. Dokl.} Vol. {\bf 32} (1985), No. {\bf 1}, 228-231.

\item[{[B]}] R. Block, On torsion-free abelian groups and Lie algebras, {\it Proc. Amer. Math. Soc.} {\bf 9} (1958), 613-620.

\item[{[CK]}] S.-J. Cheng and V. G. Kac, A new $N=6$ superconformal algebras, {\it Commun. Math. Phys.} {\bf 186} (1997), 219-231.

\item[{[F]}] V. T. Filipov, A class of simple nonassociative algebras, {\it Mat. Zametki} {\bf 45} (1989), 101-105.

\item[{[GDi1]}] I. M. Gel'fand and L. A. Dikii, Asymptotic behaviour of the resolvent of Sturm-Liouville equations and the algebra of the Korteweg-de Vries equations, {\it Russian Math. Surveys} {\bf 30:5} (1975), 77-113.

\item[{[GDi2]}] 
I. M. Gel'fand and L. A. Dikii, A Lie algebra structure in a formal variational Calculation, {\it Func. Anal. Appl.}  {\bf 10} (1976), 16-22.

\item[{[GDo1]}] 
I. M. Gel'fand and I. Ya. Dorfman, Hamiltonian operators and algebraic structures related to them, {\it Funkts. Anal. Prilozhen}  {\bf 13} (1979), 13-30.

\item[{[GDo2]}] 
I. M. Gel'fand and I. Ya. Dorfman, 
Schouten brackets and Hamiltonian operators, {\it Funkts. Anal. Prilozhen}  {\bf 14} (1980),  71-74.

\item[{[GDo3]}] 
I. M. Gel'fand and I. Ya. Dorfman, 
 Hamiltonian operators and infinite-dimensional Lie algebras, {\it Funkts. Anal. Prilozhen}  {\bf 14} (1981),  23-40.

\item[{[Ka1]}] V. G. Kac, Lie superalgebras, {\it Adv. Math.} {\bf 26} (1977), 8-96.

\item[{[Ka2]}] V. G. Kac, {\it Vertex algebras for beginners}, University lectures series, Vol {\bf 10}, AMS. Providence RI, 1996.

\item[{[Ka3]}] V. G. Kac, Superconformal algebras and transitive group actions on quadrics, {\it Commun. Math. Phys.} {\bf 186} (1997), 233-252.

\item[{[Ka4]}] V. G. Kac, Idea of locality, {\it Physical Applications and Mathematical Aspects of Geometry, Groups and Algebras}, Doebener et al eds., World Scientific Publishers, 1997, 16-32.

\item[{[KL]}]  V. G. Kac and J. W. Leur, On classification of superconformal algebras, in S. J. et al. eds. {\it String 88}, World Sci. (1989), 77-106.

\item[{[KT]}] V. G. Kac and I. T. Todorov, Superconformal current algebras and their unitary representations, {\it Commun. Math. Phys}. {\bf 102} (1985), 337-347.

\item[{[K]}] H. Kim, Complete left-invariant affine structures on nilpotent Lie
groups, {\it J. Deff. Geom.} {\bf 24} (1986), 373-394.

\item[{[M]}] P. Mathieu, Supersymmetry extension of the Korteweg-de Vries equation, {\it J. Math. Phys.} {\bf 29} (11) (1988), 2499-2507.

\item[{[O1]}]
J. Marshall Osborn, Novikov algebras, {\it Nova J. Algebra} \& {\it Geom.} {\bf 1} (1992), 1-14.

\item[{[O2]}]
J. Marshall Osborn, Simple Novikov algebras with an idempotent, {\it Commun. Algebra} {\bf 20} (1992), No. 9, 2729-2753.

\item[{[O3]}]
J. Marshall Osborn, Infinite dimensional Novikov algebras of characteristic 0, {\it J. Algebra} {\bf 167} (1994), 146-167.

\item[{[O4]}]
J. Marshall Osborn, Modules for Novikov algebras, {\it Proceeding of the II International Congress on Algebra, Barnaul, 1991.}

\item[{[O5]}]
J. Marshall Osborn, A conjecture on locally Novikov algebras, {\it Non-associative algebra and its applications} (Oviedo, 1993), 309-313, {\it Math. Appl.} {\bf 303}, Kluwer Acad. Publ., Dordrecht, 1994.

\item[{[O6]}]
J. Marshall Osborn, Modules for Novikov algebras of characteristic 0, {\it Commun. Algebra} {\bf 23} (1995), 3627-3640.

\item[{[OZ]}] J. Marshall Osborn and E. Zelmanov, Nonassociative algebras related to Hamiltonian operators in the formal calculus of variations, {\it J. Pure. Appl. Algebra} {\bf 101} (1995), 335-352.

\item[{[X1]}], X. Xu, Hamiltonian operators and associative algebras with a derivation, {\it  Lett. Math. Phys.} {\bf 33} (1995), 1-6.

\item[{[X2]}] X. Xu, Hamiltonian superoperators, {\it J. Phys A: Math. \& Gen.} {\bf 28} No. 6 (1995).

\item[{[X3]}] X. Xu, On simple Novikov algebras and their irreducible modules, {\it J. Algebra} {\bf 185} (1996), 905-934.

\item[{[X4]}] X. Xu, Novikov-Poisson algebras, {\it J. Algebra} {\bf 190} (1997), 253-279.

\item[{[X5]}] X. Xu, Skew-symmetric differential operators and combinatorial identities, {\it Mh. Math.} {\bf 127} (1999), 243-258.

\item[{[X6]}] X. Xu, Variational calculus of supervariables and related algebraic structures, {\it J. Algebra}, in press; preprint was circulated in January 1995.

\item[{[X7]}] X. Xu, {\it Introduction to Vertex Operator Superalgebras and Their Modules}, Kluwer Academic Publishers, Dordrecht/Boston/London, 1998.

\item[{[X8]}] X. Xu, Generalizations of Block algebras, {\it Manuscripta Math.}, in press.

\item[{[Z]}]
E. I. Zelmanov, On a class of local translation invariant Lie algebras, {\it Soviet Math. Dokl.} Vol {\bf 35} (1987), No. {\bf 1}, 216-218.

\end{description}
\end{document}